\documentclass[11pt]{amsart}
  \usepackage{eucal,color,extarrows,fullpage,graphics,mathptmx}
  \usepackage[shortlabels]{enumitem}
  \usepackage[all]{xy}
  \usepackage[colorlinks,citecolor=magenta,linkcolor=black]{hyperref}
 
  \SetSymbolFont{operators}{normal}{OT1}{ptm}{m}{n}  
  \DeclareSymbolFont{emcmr}{OT1}{cmr}{m}{n}   
  \DeclareMathSymbol{\Omega}{\mathalpha}{emcmr}{"0A}
  \DeclareMathSymbol{\Gamma}{\mathalpha}{emcmr}{"00}
  \DeclareMathSymbol{\Psi}{\mathalpha}{emcmr}{"09}
  \DeclareMathSymbol{\Lambda}{\mathalpha}{emcmr}{"03}

  \numberwithin{equation}{section}

  \DeclareMathOperator{\cok}{cok}
  \DeclareMathOperator{\im}{im}
  \DeclareMathOperator{\Sym}{Sym}
  \DeclareMathOperator{\Ext}{Ext}
  \DeclareMathOperator{\Hom}{Hom}
  \DeclareMathOperator{\Kos}{Kos}
  \DeclareMathOperator{\h}{\cal{H}}
  \DeclareMathOperator{\exterior}{\mbox{\large $\wedge$}}

  \newcommand{\cO}{\mathcal{O}}
  \newcommand{\isomlong}{\xlongrightarrow{\,\smash{\raisebox{-0.65ex}{\ensuremath{\scriptstyle\sim}}}\,}}
  \newcommand{\cal}[1]{\mathcal{#1}}
  \newcommand{\To}{\longrightarrow}
  \newcommand{\vsim}{\rotatebox{90}{$\displaystyle\sim$}}
  \newcommand{\id}{\mathrm{id}}

  \newcommand{\obs}{\operatorname{obs}}
  \newcommand{\class}{\operatorname{cl}}
  \newcommand{\scindage}{\operatorname{sc}}
  \newcommand{\prj}{\operatorname{pr}}
  
  \newcommand{\ctr}{\operatorname{ctr}}
  \newcommand{\product}{\mathrm{prod}}
  \newcommand{\eval}{\operatorname{ev}}
  
  \newtheorem{theorem}{Theorem}[section]
  \newtheorem{lemma}[theorem]{Lemma}
  \newtheorem{proposition}[theorem]{Proposition}
  \newtheorem{corollary}[theorem]{Corollary}

  \theoremstyle{definition} 
  \newtheorem{defin}[theorem]{Definition}
  \newtheorem{remark}[theorem]{Remark}
  \newtheorem{example}[theorem]{Example}

  \author[P.\ Achinger]{Piotr Achinger}
  \address{Institute of Mathematics of the Polish Academy of Sciences \newline \indent ul.\ Śniadeckich 8, 00-656 Warsaw, Poland}
  \email{pachinger@impan.pl}

  \author[J.\ Suh]{Junecue Suh}
  \address{1156 High Street \newline \indent University of California in Santa Cruz \newline \indent Santa Cruz CA 95064 USA}
  \email{jusuh@ucsc.edu}

  \title{Some refinements of the Deligne--Illusie theorem}
  \date{\today}

\begin{document}

\begin{abstract}
  We extend the results of Deligne and Illusie on liftings modulo $p^2$ and decompositions of the de Rham complex in several ways. We show that for a smooth scheme $X$ over a perfect field $k$ of characteristic \mbox{$p>0$}, the truncations of the de Rham complex in $\max(p-1, 2)$ consecutive degrees can be reconstructed as objects of the derived category in terms of its truncation in degrees at most one (or, equivalently, in terms the obstruction class to lifting modulo $p^2$). Consequently, these truncations are decomposable if $X$ admits a lifting to $W_2(k)$, in which case the first nonzero differential in the conjugate spectral sequence appears no earlier than on page $\max(p,3)$ (these corollaries have been recently strengthened by Drinfeld, Bhatt--Lurie, and Li--Mondal). Without assuming the existence of a lifting, we describe the gerbes of splittings of two-term truncations and the differentials on the second page of the conjugate spectral sequence, answering a question of Katz. 

  The main technical result used in the case $p>2$ belongs purely to homological algebra. It concerns certain commutative differential graded algebras whose cohomology algebra is the exterior algebra, dubbed by us \emph{abstract Koszul complexes}, of which the de Rham complex in characteristic $p$ is an example.

  In the appendix, we use the aforementioned stronger decomposition result to prove that Kodaira--Akizuki--Nakano vanishing and Hodge--de Rham degeneration both hold for $F$-split $(p+1)$-folds.
\end{abstract}

\maketitle

\section{Introduction}

\subsection{Decompositions of the de Rham complex}

In \cite{DeligneIllusie}, Deligne and Illusie showed that for a smooth scheme $X$ over a perfect field $k$ of characteristic $p>0$, a flat lifting of the Frobenius twist $X' = F_k^* X$ to $W_2(k)$ induces a splitting of the truncation of the de Rham complex in degrees $[0,1]$, i.e.\ an isomorphism in the derived category
\[ 
  \cO_{X'} \oplus \Omega^1_{X'/k}[-1]
  \isomlong 
  \tau_{\leq 1} (F_{X/k,*} \Omega^\bullet_{X/k}).
\]
Using the  algebra structure of the de Rham complex, they further show that it induces an isomorphism
\[ 
  \bigoplus_{i <p } \Omega^i_{X'/k}[-i]
  \isomlong
  \tau_{< p } (F_{X/k,*} \Omega^\bullet_{X/k}).
\]
With their method, it is unclear if one could extend this further to an isomorphism  between $\bigoplus_{i \geq 0} \Omega^i_{X'/k}[-i]$ and $F_{X/k,*} \Omega^\bullet_{X/k}$ if $\dim X\geq p$, i.e.\ whether the de Rham complex $\Omega^\bullet_{X/k}$ is \emph{decomposable}. As a step further, Deligne and Illusie prove using duality that this is the case if $\dim X = p$.

It is as of today an open problem whether there exists a smooth $X$ over $k$ liftable to $W_2(k)$, necessarily of dimension $\dim X>p$, for which the de Rham complex is not decomposable. In this paper, as a small contribution to this question, we investigate the ways in which the truncation $\tau_{\leq 1} (F_{X/k,*} \Omega^\bullet_{X/k})$ determines the truncations $\tau_{[a,b]} (F_{X/k,*}\Omega^\bullet_{X/k})$. Our first result is the following:

\begin{theorem} \label{thm:main-intro}
  Let $X$ be a smooth scheme over a perfect field $k$ of characteristic $p>0$ which is liftable to $W_2(k)$. Then the truncations
  \[ 
    \tau_{[a,b]} (F_{X/k,*}\Omega^\bullet_{X/k})
  \]
  are decomposable for $a\leq b < a+p-1$ when $p>2$, and for $a\leq b \leq a+1$ when $p=2$.
\end{theorem}

The above result immediately implies that in the conjugate spectral sequence
\begin{equation}  \label{eqn:conj-ss-intro}
  E_2^{ij} = H^i(X', \Omega^j_{X'/k}) 
  \quad \Rightarrow \quad 
  H^{i+j}(X, \Omega^\bullet_{X/k})
\end{equation}
the differentials $d_r^{ij}$ are zero for $2\leq r<p$ when $p>2$, and for $r=2$ when $p=2$. As a sample corollary, we obtain the following criterion for degeneration of spectral sequences.

\begin{corollary} \label{cor:ss-deg-intro}
  For $X$ as in Theorem~\ref{thm:main-intro}, suppose that 
  \[ 
    H^i(X, \Omega^j_{X/k}) = 0
    \quad \text{for } \quad
     |i-j| \geq p. 
  \]
  Then the conjugate spectral sequence \eqref{eqn:conj-ss-intro} degenerates. If moreover $X$ is proper over $k$, then the Hodge to de Rham spectral sequence
  \[ 
    E^{ij}_1 = H^j(X, \Omega^i_{X/k})
    \quad \Rightarrow \quad
    H^{i+j}(X, \Omega^\bullet_{X/k})
  \]
  degenerates as well. 
\end{corollary}

\subsection{Truncations of the de Rham complex}

Our methods give information about truncations of the de Rham complex without assuming liftability modulo $p^2$. Our results in this direction are the strongest and most explicit for truncations in two consecutive degrees. Namely, for a general smooth $X$ over $k$ (not necessarily liftable to $W_2 (k)$) and for $q\geq 1$, the truncated complex $\tau_{[q-1, q]}(F_{X/k,*}\Omega^\bullet_{X/k})$ can be described as the mapping fiber of $\delta^q[-q]$ for a map
\[ 
  \delta^q\colon \Omega^q_{X'/k} \To \Omega^{q-1}_{X'/k}[2]
  \quad \text{i.e.} \quad
  \delta^q \in \Ext^2(\Omega^q_{X'/k}, \Omega^{q-1}_{X'/k})
\]
which is the ``cup product'' with the negative of the deformation obstruction class
\begin{equation} \label{eqn:def-class}
  \delta^1 = -\obs(X'/k/W_2(k)) \quad \in \quad \Ext^2(\Omega^1_{X'/k}, \cO_{X'}) \simeq H^2(X', T_{X'/k})
\end{equation}
to the existence of a lifting of $X'$ to $W_2(k)$ (see Corollary~\ref{cor:delta-q-formula}). The result in particular implies that the two-term truncation $\tau_{[q-1,q]}(F_{X/k,*}\Omega^\bullet_{X/k})$ is decomposable if $X'$ lifts to $W_2(k)$, and yields a description the differentials on the second page of the conjugate spectral sequence --- answering a natural question of Katz.

\begin{theorem}[see Corollary~\ref{cor:katz-question}] \label{thm:katz-conj-intro}
  In the above situation, the differential
  \[
    d_2^{ij} \colon H^i(X', \Omega^j_{X'/k}) \To H^{i+2}(X', \Omega^{j-1}_{X'/k})
  \]
  in the conjugate spectral sequence \eqref{eqn:conj-ss-intro} is induced by the cup product with the negative of the obstruction class $\obs(X'/k/W_2(k))$.
\end{theorem}

Deligne and Illusie \cite[\S 3]{DeligneIllusie} define the \emph{gerbe of splittings} $\scindage K$ of a two-term complex $K$, and relate the gerbe of splittings of $\tau_{\leq 1}(F_{X/k,*}\Omega^\bullet_{X/k})$ to the gerbe of liftings of $X'$ to $W_2(k)$. This provides a ``categorification'' of the equality \eqref{eqn:def-class}. In the same vein, for $p>2$, our description of the class of $\tau_{[q-1,q]}(F_{X/k,*}\Omega^\bullet_{X/k})$ can be upgraded to a morphism of gerbes (see Theorem~\ref{thm-s2.2.1})
\[ 
  \wedge^q \colon \scindage(\tau_{\leq 1}(F_{X/k,*}\Omega^\bullet_{X/k})) \To \scindage(\tau_{[q-1,q]}(F_* \Omega^\bullet_{X/k})).
\]

Let us now discuss longer truncations of the de Rham complex. The assertion of Theorem~\ref{thm:main-intro} is subsumed by a recent beautiful observation of Drinfeld \cite[\S 5.12.1]{Drinfeld} (a proof will appear in Bhatt and Lurie's forthcoming \cite{BhattLurie}): a lifting of $X'$ to $W_2(k)$ induces a $\mu_p$-action on the de Rham complex $F_{X/k,*} \Omega^\bullet_{X/k}$, which one can use to show that the truncations $\tau_{[q-p+1, q]}(F_{X/k,*} \Omega^\bullet_{X/k})$ are decomposable for all $q$ (even more recently, Li and Mondal \cite{LiMondal} found an independent proof). However, the method of proof of Theorem~\ref{thm:main-intro} is completely different and provides interesting information even if $X$ is not liftable to $W_2(k)$. It is deduced from the following result (when $p>2$) and Corollary \ref{cor:delta-q-formula} (when $p=2$) alluded to above.

\begin{theorem}\label{thm:main-intro-2}
  Let $X$ be a smooth scheme over a perfect field $k$ of characteristic $p>0$, let $q$ be an integer, and let $m<p$. One then has an isomorphism in the derived category of $X'$:
  \[
    \tau_{[q-m, q]}(F_{X/k, *}\Omega^\bullet_{X/k}) \quad\simeq\quad \tau_{\geq q-m}(L\Gamma^q(\tau_{\leq 1} F_{X/k, *}\Omega^\bullet_{X/k}))
  \]
  where $L\Gamma^q$ is the derived $q$-th divided power.
\end{theorem}

\subsection{Abstract Koszul complexes}

The proof of Theorem~\ref{thm:main-intro-2} has very little to do with algebraic geometry. To state the main technical result behind it, we need the notion of an \emph{abstract Koszul complex} (Definition~\ref{def:abstr-koszul-cx}), which is a certain commutative differential graded algebra (cdga) $K$ in a ringed topos for which the multiplication induces isomorphisms
\[
  \exterior^q \h^1(K) \isomlong \h^q(K) 
  \quad \text{for all} \quad
  q\geq 0.
\]
Thanks to the Cartier isomorphism, the de Rham complexes $F_{X/k,*}\Omega^\bullet_{X/k}$ in characteristic $p>0$ are examples of such, and hence the result below immediately implies Theorem~\ref{thm:main-intro-2}.

\begin{theorem}[see Theorem~\ref{thm:koszul-reconstr}] \label{thm:koszul-intro}
  Let $K$ be an abstract Koszul complex in a ringed topos $(X, \cO)$ satisfying the flatness condition \eqref{eqn:flatness}, and let $q\geq m\geq 1$ be integers such that $m!$ is invertible in $\cO$. Suppose that either $q=m$ or that $m+1$ is a~nonzerodivizor in $\cO$. Then there exists an isomorphism in the derived category
  \begin{equation} \label{eqn:koszulmap-intro}
    \tau_{\geq q-m} \left( L\Gamma^q (\tau_{\leq 1}(K)) \right) \isomlong \tau_{[q-m, q]}(K).
  \end{equation}
\end{theorem}

\noindent 
In \eqref{eqn:koszulmap-intro}, $L\Gamma^q$ is again the derived $q$-th divided power, and the source of the map can be more concretely realized as $\tau_{\geq q-m}$ of the Koszul complex 
\[ 
    \cdots \to \underbrace{\exterior^i (\cal{Z}^1 K) \otimes \Gamma^{q-i}(K^0)}_{\text{degree $i$}} \to \cdots
    \to \exterior^{q-1}(\cal{Z}^1 K) \otimes K^0 \to \underbrace{\exterior^q (\cal{Z}^1 K)}_{\text{degree }q} \to 0 \to \cdots  
\]

For $m=1$ (and assuming that $2$ is a nonzerodivisor), we again obtain more refined information regarding $\tau_{[q-1, q]} K$, including the differentials on the second page of the spectral sequence
\[ 
    E_2^{ij} = H^i(X, \h^j(K)) \quad \Rightarrow \quad H^{i+j}(X, K)
    \qquad \text{(see Corollary~\ref{cor:delta-q-formula})}.
\]

As observed by Kato \cite{Kato}, logarithmic de Rham complexes are abstract Koszul complexes, and hence Theorem~\ref{thm:main-intro} works also in the log case. The inspiration for Theorem~\ref{thm:koszul-reconstr} came from the result of Steenbrink \cite[\S 2.8]{Steenbrink} describing the nearby cycle complex $R\Psi \mathbf{Q}$ for a complex semistable degeneration in terms of the logarithmic structure; see also \cite[\S 4]{AchingerOgus}. It is an interesting question whether Steenbrink's result can be extended to work with integral coefficients; the nearby cycles $R\Psi \mathbf{Z}$ are co-connective $E_\infty$-algebra versions of abstract Koszul complexes, but we do not know whether they admit cdga models (see Remark~\ref{rmk:einfty} and Example~\ref{ex:nearby}). An affirmative answer would give an application unrelated to the Deligne--Illusie theorem, refining \cite[Theorem~4.2.2(1)]{AchingerOgus}, providing a description of the two-step truncations $\tau_{[q-1,q]}$ of certain logarithmic nearby cycle complexes.

\subsection{\texorpdfstring{The case $p=2$ (Theorem \ref{thm-s2.3.1})}{The case p=2}}

The description of the truncations $\tau_{[q-1,q]}(F_{X/k,*} \Omega^\bullet_{X/k})$ and its corollary, Theorem~\ref{thm:katz-conj-intro}, can be deduced from the ``abstract Koszul complex'' machinery and Theorem~\ref{thm:main-intro-2}, but only for $p>2$. In contrast, the assertion of Theorem~\ref{thm:koszul-intro} is vacuous if $2\cdot \cO=0$. Accordingly, the computation of the class of $\tau_{[q-1,q]} F_{X/k,*} \Omega^\bullet_{X/k}$, occupying the entire Section~\ref{s:gerbe-calculation} is much harder in the case $p=2$, and uses more information about the de Rham complex than merely its abstract Koszul complex structure.  For this technical point, we highlight the passage from (\ref{eqn-s2.7}) to (\ref{eqn-s2.7bis}).

It could be worthwile to extend the methods used in the case of $p=2$ in order to ``compute'' the truncations $\tau_{[q-p+1, p]}$ in $p$ consecutive degrees, and it would be interesting to extract the exact abstract properties of the de Rham complex in positive characteristic needed for the proof. Its relationship with the aforementioned result of Drinfeld, Bhatt--Lurie, and Li--Mondal remains elusive.

The results concerning the truncations $\tau_{[q-1,q]} (F_{X/k,*} \Omega^\bullet_{X/k})$ and Theorem~\ref{thm:katz-conj-intro}, including the case $p=2$, presented here are due to the second author and appeared in his 2007 Ph.D.\ thesis. After the first author proved Theorem~\ref{thm:main-intro-2}, the authors decided to publish their results together.  

\subsection{\texorpdfstring{Application to $F$-split $(p+1)$-folds}{Application to F-split (p+1)-folds}}

As an illustration of this circle of ideas, using the refinement of the Deligne--Illusie theorem due to Drinfeld, Bhatt--Lurie, and Li--Mondal, we prove in Appendix~A that the Kodaira--Akizuki--Nakano vanishing theorem and the degeneration of the Hodge to de Rham spectral sequence both hold for $F$-split $(p+1)$-folds in characteristic $p$.

\subsection*{Acknowledgements} 

The first author is grateful to Luc Illusie, Shizhang Li, Arthur Ogus, and Vadim Vologodsky for useful comments and stimulating discussions. 
The second author thanks Nicholas Katz for raising the question that we answer in Theorem \ref{thm:katz-conj-intro},  and Pierre Deligne and Luc Illusie for helpful comments.

The first author was supported by NCN OPUS grant number {\tt UMO-2015/17/B/ST1/02634}.

\subsection*{Notation}

If $K$ is a cochain complex in an abelian category, we write $\cal{Z}^i K = \ker(d\colon K^i\to K^{i+1})$, $\cal{B}^i K = \im(d\colon K^{i-1}\to K^i)$, and $\h^i(K)=\cal{Z}^i K/\cal{B}^i K$, and denote by $\tau_{\leq b} (K)$ the subcomplex
\[ 
  \tau_{\leq b}(K) = \left[ \cdots \to K^{b-1} \to \cal{Z}^b K \to 0 \to \cdots \right],
\]
by $\tau_{\geq a}(K)$ the quotient $K/\tau_{\leq a}(K)$, and define $\tau_{[a,b]}(K) = \tau_{\geq a} \tau_{\leq b} (K)$. We call $K$ \emph{decomposable} if it is isomorphic in the derived category to the complex with zero differential $\bigoplus \h^i(K)[-i]$. 

A \emph{commutative differential graded algebra} (\emph{cdga}) is an associative graded ring $K = \bigoplus_{n\in \mathbf{Z}} K^n$ which is graded-commutative (i.e.\ $xy = (-1)^{mn} yx$ for $x\in K^m$, ${y\in K^n}$), endowed with a differential $d\colon K\to K$ mapping $K^n$ to $K^{n+1}$ and satisfying $d(xy) = dx\cdot y + (-1)^n x \cdot dy$ for $x\in K^n$. We say that $K$ is \emph{coconnective} if $K^n = 0$ for $n<0$.

\section{Abstract Koszul complexes}

\subsection{Definition and examples} We work in a ringed topos $(X, \cO)$. 

\begin{defin}[Abstract Koszul complex] \label{def:abstr-koszul-cx}
  A coconnective commutative differential graded $\cO$-algebra $K$ is called an \emph{abstract Koszul complex} if the following conditions are satisfied:
  \begin{enumerate}[(i)]
    \item \emph{$\cO \to \h^0(K)$ is an isomorphism,}
    \item \emph{For every $q\geq 1$, the multiplication map $\h^1(K)^{\otimes q}\to \h^q(K)$ factors through an isomorphism}
      \[ 
        \mu^q : \exterior^q \h^1(K) \isomlong \h^q(K).
      \]
  \end{enumerate}
\end{defin}

\begin{example}[De Rham complex in characteristic $p>0$] \label{ex:derham-koszul}
  Let $X$ be a smooth scheme over a perfect field $k$ of characteristic $p>0$, and let $F_{X/k}\colon X\to X'$ be its relative Frobenius. Let $K = F_{X/k,*}\Omega^\bullet_{X/k}$ be the de Rham complex, treated as a cdga over $\cO_{X'}$. Then the Cartier isomorphisms 
  \[ 
    C\colon \h^i(F_{X/k,*}\Omega^\bullet_{X/k}) \isomlong \Omega^i_{X'/k}  
  \]
  are multiplicative, and hence $K$ is an abstract Koszul complex over $(X', \cO_{X'})$. 

  More generally, if $f\colon (X, \cal{M}_X)\to (S, \cal{M}_S)$ is a morphism of fine log schemes over $\mathbf{F}_p$ which is smooth and of Cartier type, then the log de Rham complex $F_{X/S,*}\Omega^\bullet_{(X, \cal{M}_X)/(S, \cal{M}_S)}$ is an abstract Koszul complex \cite[Theorem~4.12]{Kato}.
\end{example}

\begin{example}[{Nearby cycle complexes, see e.g.\ \cite[\S 2]{Steenbrink}}] \label{ex:nearby}
  Let $X$ be a complex manifold and let $D = \bigcup D_\alpha$ be a divisor with simple normal crossings on $X$. Let $j\colon U=X\setminus D \hookrightarrow X$ be the complementary open immersion, and let $K = Rj_* \mathbf{Q}_U$. Since we are working with rational coefficients, we can find a cdga model for $K$ (e.g.\ \cite[Part~II, Corollary~1.5]{KrizMay}). The purity theorem implies that
  \[ 
    R^i j_* \mathbf{Q}_U = \exterior^i \left( \bigoplus \mathbf{Q}_{D_\alpha} \right), 
  \]
  and hence any cdga model of $K$ is an abstract Koszul complex over $(X, \mathbf{Q}_X)$. Moreover, one has an isomorphism in the derived category \cite[Lemma~2.7]{Steenbrink} (see also \cite[\S 4]{AchingerOgus})
  \[ 
    \tau_{\leq 1} (R^1 j_* \mathbf{Z}_U) 
    \quad \simeq \quad
    \left[ \cO_X \xrightarrow{\exp} \cal{M}^{\rm gp} \right],
  \]
  where $\cal{M}^{\rm gp}$ is the sheaf of meromorphic functions without zeros or poles on $U$. Variants of this contruction exist for the nearby cycle complexes $R\Psi \mathbf{Q}$ for a semistable degeneration over a disc, and there exist analogs in $\ell$-adic \'etale cohomology.
\end{example}

In the following, we will work with abstract Koszul complexes satisfying the additional flatness condition:
\begin{equation} \label{eqn:flatness}
  \text{\emph{The $\cO$-modules $K^0$, $\cal{B}^1 K$, $\cal{Z}^1 K$, and $\h^1(K)$ are flat.}}
\end{equation}
In particular, this implies that the modules $\h^q(K) \simeq \exterior^q \h^q(K)$ are flat for all $q\geq 0$. The above condition is satisfied in the situation of Examples~\ref{ex:derham-koszul} and \ref{ex:nearby}. 

\subsection{Koszul complexes}
\label{section-Koszul-complexes}

Our goal is to show that to a certain extent, the underlying complex of an abstract Koszul complex satisfying the flatness condition \eqref{eqn:flatness} is determined by its truncation in degrees $\leq 1$ (Theorem~\ref{thm:koszul-reconstr}). We achieve this using the notion of the Koszul complex of a~map $u:P\to Q$, see \cite[Chapitre I, \S 4.3]{Illusie} and \cite[\S 1.1--1.2]{KatoSaito}.

Recall first that if $M$ and $N$ are $\cO$-modules over, then for every $q\geq 0$ there is a natural decomposition of the divided (resp.\ exterior) power 
\[
  \Gamma^q (M\oplus N) = \bigoplus_{a+b=q} \Gamma^a M \otimes \Gamma^b N
  \qquad
  \text{(resp.\ $\exterior^q (M\oplus N) = \bigoplus_{a+b=q} \exterior^a M \otimes \exterior^b N$)}
\]
In what follows, we will use the \emph{comultiplication maps} 
\[
  \eta^q:\Gamma^q M \To M\otimes \Gamma^{q-1} M
  \qquad
  \text{(resp.\ $\eta^q:\exterior^q M \To (\exterior^{q-1} M)\otimes M)$} 
\]
obtained as the composition of $\Gamma^q$ (resp.\ $\exterior^q$) of the diagonal map $M\to M\oplus M$ and the projection to the $(a, b) = (1, q-1)$-part (resp.\ $(a, b)=(q-1, 1)$-part) in the above decomposition of $\Gamma^q(M\oplus M)$ (resp.\ of $\exterior^q(M\oplus M)$). Explicitly, we have
\begin{align*}
  \eta^q(x_1^{[e_1]}\cdot \ldots \cdot x_r^{[e_r]}) &= \sum_{i=1}^r x_i \otimes (x_1^{[e_1]} \cdot\ldots \cdot x_i^{[e_i - 1]} \cdot \ldots \cdot x_r^{[e_r]}), & (e_1 + \cdots + e_r = q) \\
  \eta^q(x_1\wedge\ldots\wedge x_q) &= \sum_{i=1}^q (-1)^{i-1} x_1\wedge \ldots \wedge \widehat{x_i}\wedge\ldots\wedge x_q \otimes x_i.
\end{align*}
Sometimes we omit the superscript $q$ when it is clear from the context.

\begin{defin}[Koszul complex $\Kos^q(u)$]
  Let $u:P\to Q$ be a map of $\cO$-modules, and let $q\geq 0$ be an integer. Then the \emph{$q$-th Koszul complex} $\Kos^q(u)$ is the cochain complex whose $i$-th term is 
  \[
    \Kos^q(u)^i = \exterior^i(Q) \otimes \Gamma^{q-i}(P), 
  \]
  with differential $d\colon \Kos^q(u)^i \to \Kos^q(u)^{i+1}$ defined as the composition 
  \[ 
    \xymatrix{
      \exterior^i(Q)\otimes \Gamma^{q-i}(P)\ar[r]^-{1\otimes \eta} &  \exterior^i(Q)\otimes P\otimes \Gamma^{q-i-1}(P) \ar[d]^{1\otimes u \otimes 1} & \\
      & \exterior^i(Q)\otimes Q \otimes \Gamma^{q-i-1}(P) \ar[r]_{\wedge \otimes 1} & \exterior^{i+1}(Q)\otimes \Gamma^{q-i-1}(P). 
    }
  \]
\end{defin}

\noindent 
Concretely, with $e_1 + \cdots + e_r = q-i$, $x_1, \ldots, x_r\in P$, and $y\in \exterior^i Q$:
\[ 
  d\left(y \otimes x_1^{[e_1]} \cdot\ldots \cdot x_r^{[e_r]}\right) 
  = 
  \sum_{j=1}^r y \wedge u(x_j) \otimes 
  \left( x_1^{[e_1]} \cdot\ldots \cdot x_j^{[e_j-1]} \cdot \ldots \cdot x_r^{[e_r]} \right).
\]

\begin{proposition}
  Let $u\colon P\to Q$ be a map of flat $\cO$-modules, and let $F(u) = [P\xrightarrow{u} Q]$ be the two-term cochain complex with $P$ in degree zero (the mapping fiber). There exist natural isomorphisms in the derived category
  \[ 
    \Kos^q(u) 
    \quad \simeq \quad
    L\Lambda^q(F[1])[-q] 
    \quad \simeq \quad 
    L\Gamma^q(F)
  \]
  where $L\Lambda^q$ (resp.\ $L\Gamma^q$) is the derived exterior (resp.\ divided) power.
\end{proposition}

\begin{proof}
Combine \cite[Corollary~1.2.7]{KatoSaito} with \cite[I 4.3.2.1]{Illusie}.
\end{proof}

\begin{corollary}
  For a map $u\colon P\to Q$ between flat $\cO$-modules, the complex $\Kos^q(u)$, treated as an object of the derived category, depends only on $F(u) = [P\to Q]$ up to quasi-isomorphism. In particular, if $F(u)$ is decomposable, then so is $\Kos^q(u)$.
\end{corollary}

\begin{proposition}[{cf.\ \cite[Lemma~1.4]{Steenbrink}}] \label{prop:alphamaps}
  Let $u\colon P\to Q$ be a map of $\cO$-modules. There exist unique arrows 
  \[
    \exterior^i(Q)\otimes\Gamma^{q-i}(\ker u) \To \cal{Z}^i(\Kos^q(u))
    \quad \text{and} \quad
     \alpha^{i} : \exterior^i(\cok u)\otimes \Gamma^{q-i}(\ker u) \To \h^i(\Kos^q(u))
  \]
  making the following diagram commute
  \[ 
    \xymatrix{
      \exterior^i(Q) \otimes \Gamma^{q-i}(P) \ar@{=}[r] & \Kos^q(u)^i \\
      \exterior^i(Q)\otimes\Gamma^{q-i}(\ker u) \ar[r] \ar[d] \ar[u] & \cal{Z}^i(\Kos^q(u)) \ar[u] \ar[d] \\
      \exterior^i(\cok u)\otimes \Gamma^{q-i}(\ker u) \ar[r]_-{\alpha^{i}} & \h^i(\Kos^q(u)).
    }
  \]
  Moreover, the map $\alpha^i$ is an isomorphism if $P$, $Q$, $\ker u$, $\im u$, and $\cok u$ are all flat.
\end{proposition}

\begin{proof}
The first assertion is straightforward. The second is reduced as in \cite[I~4.3.1.6]{Illusie} to the case where $P$, $Q$, $\ker u$, $\im u$, and $\cok u$ are free $\cO$-modules of finite rank. In this case, splitting the surjection $Q\to \cok u$ one can write $u = u'\oplus u''$ where $u'\colon P\to \im u$ and $u''\colon 0\to \cok u$. The assertion then holds for $u'$ (by \cite[4.3.1.6]{Illusie}) and for $u''$ (trivially), for all $q$, and then the assertion for $u = u'\oplus u''$ follows from the isomorphism \cite[4.3.1.5]{Illusie} 
\[ 
  \Kos^\bullet(u) = \Kos^\bullet(u') \otimes \Kos^\bullet(u'')
\]
where $\Kos^\bullet(u) = \bigoplus_{q\geq 0} \Kos^q(u)[q]$.
\end{proof}

\subsection{Truncations of abstract Koszul complexes}

The following theorem is the main result of this section.

\begin{theorem} \label{thm:koszul-reconstr}
  Let $m$ be an integer such that $m!$ is invertible in $\cO$, and let $q\geq m$. Suppose that either $q= m$, or that $m+1$ is not a zero divisor in $\cO$.   Let $K$ be an abstract Koszul complex on $(X, \cO)$ satisfying the flatness condition \eqref{eqn:flatness}, and write 
  \[
    \tau_{\leq 1} K = \left[ K^0 \xrightarrow{\partial} Z^1 K \right] 
  \] 
  for its truncation in degrees $\leq 1$. Then the multiplication maps
  \[ 
    \Kos^q(\partial)^i =  \exterior^i(\cal{Z}^1 K) \otimes \Gamma^{q-i}(K^0) = \exterior^i(\cal{Z}^1 K) \otimes \Sym^{q-i}(K^0) \To K^{i} 
  \]
  for $q-m\leq i \leq q$ (where we can identify $\Gamma^{q-i}$ with $\Sym^{q-i}$ as $q-i\leq m$, so that $(q-i)!$ is invertible in $\cO$) induce a quasi-isomorphism
  \begin{equation} \label{eqn:reconstr}
    \tau_{\geq q-m}(L\Gamma^q(\tau_{\leq 1}(K)) = \tau_{\geq q-m} \Kos^q(\partial) \isomlong \tau_{[q-m, q]}(K).
  \end{equation}
\end{theorem}

\begin{proof}
The multiplication maps define a morphism of ``naive truncations'' 
\[ 
  \xymatrix@C=1.2em{
    \Kos^q(\partial)^{\geq q-m}=[ \ar[d]_\mu  & \exterior^{q-m}(\cal{Z}^1 K) \otimes \Sym^m (K^0) \ar[d]\ar[r] & \cdots \ar[r] & \exterior^q (\cal{Z}^1 K)  \ar[d] & ] \\
    \tau_{\leq q}(K)^{\geq q-m}=[  & K^{q-m} \ar[r] & \cdots \ar[r] & Z^q K & ]. 
  }
\]
To obtain the desired morphism $\mu\colon \tau_{\geq q-m} \Kos^q(\partial) \to \tau_{[q-m,q]}(K)$, we need to check that the map 
\[
  \Kos^q(\partial)^{q-m} \To K^{q-m}
\] 
takes the image of $\Kos^q(\partial)^{q-m-1}\to \Kos^q(\partial)^{q-m}$ into $\cal{B}^{q-m} K = dK^{q-m-1}$. This is clear if $q=m$, so suppose that $(m+1)$ is not a zero divisor. 

Let $z\in \Kos^q(\partial)^{q-m}$ be the image of $w\in \Kos^q(\partial)^{q-m-1}$, and consider $(m+1)w$ as an element of the submodule
\[
  \exterior^{q-m-1}(\cal{Z}^1 K) \otimes \Sym^{m+1}(K^0)
  \quad \subseteq \quad
  \exterior^{q-m-1}(\cal{Z}^1 K) \otimes \Gamma^{m+1}(K^0).
\] 
Let $u\in K^{q-m-1}$ be the image of $(m+1)w$ under the multiplication map 
\[
  \exterior^{q-m-1}(\cal{Z}^1 K)\otimes \Sym^{m+1}(K^0)  \To K^{q-m-1}.
\] 
Then $du = (m+1)\mu(z)$ in $K^{q-m}$, where $\mu(z)$ is the image of $z$ under the multiplication map, and hence $\mu(z)$ gives an $(m+1)$-torsion class in $\h^{q-m}(K)$. Since by assumption $\h^{q-m}(K)\simeq \exterior^{q-m} \h^1(K)$ is flat and $m+1$ is not a zero divisor, $\mu(z)\in dK^{m-q-1}$ as desired. 

Finally, the maps induced by $\mu\colon \tau_{\geq q-m} \Kos^q(\partial) \to \tau_{[q-m,q]}(K)$ on cohomology can, thanks to Proposition~\ref{prop:alphamaps}, be identified with the maps 
\[ 
  \mu^i \colon \exterior^i \h^1(K) \To \h^i(K) \quad \text{for} \quad q-m\leq i \leq q, 
  \]
which are isomorphisms by assumption.
\end{proof}

\begin{remark}
Implicit in the above proof is the subcomplex $\widetilde{\Kos^q}(u)$ of $\Kos^q(u)$ whose $i$-th term equals $\exterior^i (Z^1 K) \otimes \Sym^{q-i} K$. The two complexes agree in degrees $\geq q-m$, and more generally the quotient $\Kos^q(u)^{q-i}/\widetilde{\Kos^q}(u)^{q-i}$ is annihilated by $i!$. This subcomplex probably does not have any ``derived meaning,'' (for example, it is not clear that it is decomposable if $\tau_{\leq 1}(K)$ is), but its advantage is that there is a~multiplication map $\mu\colon \widetilde{\Kos^q}(u)\to K$.
\end{remark}

\begin{remark} \label{rmk:einfty}
Our proof of Theorem~\ref{thm:koszul-reconstr} makes use of an explicit model of the cdga $K$. Thus, for example, if $K$ and $K'$ are equivalent cdgas to which the theorem applies, it is not obvious whether the isomorphisms \eqref{eqn:reconstr} we obtain for $K$ and $K'$ are compatible. More importantly, it does not apply to the more general case of coconnective $E_\infty$-algebras or cosimplicial commutative rings whose cohomology algebras satisfy axioms (i)--(ii) of Definition~\ref{def:abstr-koszul-cx}.
\end{remark}

\begin{corollary} \label{cor:abstract-decomp}
  Let $K$ be an abstract Koszul complex, and let $m$ be such that $m!$ is invertible in $\cO$. Suppose that $\tau_{\leq 1}(K)$ is decomposable. Then  for $a\leq b < a+m$, the complex $\tau_{[a,b]}(K)$ is decomposable. Moreover, the complex $\tau_{\leq m}(K)$ is decomposable as well.
\end{corollary}

\subsection{Application to de Rham cohomology}

We now establish some of the straightforward consequences for de Rham cohomology mentioned in the introduction. The remaining ones shall be established at the end of Section~\ref{s:gerbe-calculation}. 

\begin{proof}[Proof of Theorem~\ref{thm:main-intro}]
Let $K = F_{X/k, *} \Omega^\bullet_{X/k}$. By Example~\ref{ex:derham-koszul}, this is an abstract Koszul complex over the ringed space $(X', \cO_{X'})$. By \cite[Th\'eor\`eme~2.1]{DeligneIllusie}, the liftability assumption implies that the complex $\tau_{\leq 1}(K)$ is decomposable. Corollary~\ref{cor:abstract-decomp} with $m=p-1$ implies that $\tau_{[a,b]}(K)$ is decomposable for $a \leq b < a+ p -1$, as desired. 
\end{proof}

\begin{proof}[Proof of Corollary~\ref{cor:ss-deg-intro}]
The differentials on the $E_r$-page of \eqref{eqn:conj-ss-intro} depend only on the truncations $\tau_{[a,b]}(\Omega^\bullet_{X/k})$ with $a\leq b < a+r$, and hence all differentials on the pages $E_r$ with $r<p$ vanish. Suppose that $d^{ij}_r\neq 0$, then in particular $H^i(X, \Omega^j_{X/k})$ and $H^{i+r}(X, \Omega^{j-r+1}_{X/k})$ are both nonzero, and hence 
\[
  |i-j| < p 
  \quad \text{and} \quad
  |(i+r)-(j-r+1)| = |(i-j) + 2r - 1| < p,
\]
which implies $r<p$. Therefore \eqref{eqn:conj-ss-intro} degenerates. 

For $X$ proper over $k$, one can deduce the degeneration of the Hodge to de Rham spectral sequence as in \cite[Corollaire~2.4]{DeligneIllusie}.
\end{proof}

\begin{remark}
    As in \cite[\S 4]{DeligneIllusie} and \cite[Theorem~4.12(2)]{Kato}, analogous assertions hold for a smooth and separated morphism of $\mathbf{F}_p$-schemes $X\to S$, or more generally for a smooth morphism of Cartier type $f\colon (X, \cal{M}_X)\to (S, \cal{M}_S)$ between fine log schemes over $\mathbf{F}_p$, assuming that there exists a fine log scheme $(\widetilde{S}, \cal{M}_{\widetilde{S}})$ over $\mathbf{Z}/p^2\mathbf{Z}$ such that $\widetilde{S}$ is flat over $\mathbf{Z}/p^2\mathbf{Z}$ and a smooth lifting 
    \[ 
        \tilde f'\colon (\widetilde X{}', \cal{M}_{\widetilde X{}'})\to (\widetilde S, \cal{M}_{\widetilde S})
    \]
    of $f'$ (the base change $f$ under the absolute Frobenius $(S, \cal{M}_S)\to (S, \cal{M}_S)$). Here $\mathbf{F}_p$ and $\mathbf{Z}/p^2\mathbf{Z}$ are given the trivial log structure. 
\end{remark}

\section{Truncations in two consecutive degrees and gerbes of splittings}

In the following, we make a more detailed analysis of the truncations $\tau_{[q-1,q]} K$ for an abstract Koszul complex $K$, as well as their associated gerbes of splittings. We keep working in the category of modules in a ringed topos $(X, \mathcal{O})$.

\subsection{First order attachment maps}

For a complex $K$ and an integer $q$, the truncation
\[ 
  \tau_{[q-1,q]} K = [ \cdots \to 0 \to K^{q-1}/\cal{B}^{q-1} \to \cal{Z}^q K \to 0\to \cdots]
\]
fits inside the functorial exact triangle
\[ 
  \h^{q-1}(K)[1-q]\To \tau_{[q-1,q]} K\To \h^q(K)[-q] \xrightarrow{\delta^q_K[-q]} \h^{q-1}(K)[2-q]
\]
yielding a morphism 
\[
  \delta^q_K \colon \h^q(K)\To \h^{q-1}(K)[2].
\]
such that $\delta^q_K[-q]$ is the unique morphism making the above triangle distinguished (see \cite[Proposition~2.1.1]{AchingerOgus}). Thus the truncation $\tau_{[q-1,q]} K$ is determined by the map $\delta^q_K$, as the mapping fiber of $\delta^q_K[-q]$; it is decomposable if and only if $\delta^q_K=0$. We note for future reference the effect of the shift functor on the maps~$\delta^q_K$:
\begin{equation} \label{eqn:delta-shift-sign}
  \delta^q_{K[p]} = (-1)^p\delta^{p+q}_K.
\end{equation}

The maps $\delta^q_K$ describe the differentials on the second page of the spectral sequence
\[ 
  E_2^{pq} = H^p(X, \h^q(K)) \quad \Rightarrow \quad H^{p+q}(X, K).
\]
Namely, the differential 
\[
  d_2^{pq}\colon H^p(X, \h^q(K)) \To H^{p+2}(X, \h^{q-1}(K)) = H^p(X, \h^{q-1}(K)[2]) 
\]
is the map induced by $\delta^q_K$ on $H^p(X, -)$.

\subsection{Gerbe of splittings}
\label{ss:gerbe-of-splittings}

We recall the gerbe of splittings described in \cite{DeligneIllusie}. Let 
\[
  K=[K^0\xrightarrow{d} K^1]
\]
be a two-term complex (i.e.\ $K^i=0$ for $i\neq 0,1$), and suppose that the two conditions below hold
\begin{enumerate}
  \item ${\cal{H}}^1 (K)$ is locally free of finite rank, and
  \item the projection of $K^0$ onto ${\cal{B}}^1 = \im (d)$ locally admits a section.
\end{enumerate}
One then constructs the gerbe $\scindage(K)$ under $\underline{\Hom}(\h^1(K), \h^0(K))$ over $X$ \cite[\S 3.2]{DeligneIllusie} as the stackification of the prestack $\scindage'(K)$ whose objects are local splittings 
\[
  s\colon \h^1(K)\To K^1
\]
of the projection $K^1=\cal{Z}^1 K\to \h^1(K)$, and where morphisms $s\to s'$ are maps 
\[
  h\colon \h^1(K)\To K^0
\]
such that $dh = s'-s$. The automorphisms of an object $s$ are then identified with $\underline{\Hom}(\h^1(K), \h^0(K))$, and this makes $\scindage(K)$ into a gerbe under $\underline{\Hom}(\h^1(K), \h^0(K))$. We denote by
\[ 
  \class \scindage K \quad \in \quad H^2(X, \underline{\Hom}(\h^1(K), \h^0(K))) = \Ext^2(\h^1(K), \h^0(K)) = \Hom(\h^1(K), \h^0(K)[2])
\]
the class of the gerbe $\scindage K$. The following result relates this class to the map $\delta^1_K$ defined previously.

\begin{lemma}[{\cite[Proposition~3.3]{DeligneIllusie}}] \label{lemma:di-scindage-class}
  Let $K=[K^0\to K^1]$ be a two-term complex satisfying (1) and (2) above. Then, one has the following equality in $\Hom(\h^1(K), \h^0(K)[2])$:
  \[
    \class \scindage K = -\delta^1_K.
  \] 
\end{lemma}

A bit more generally, given an integer $q$ and a complex $K$ satisfying the following conditions:
\begin{enumerate}
  \item $K^i = 0$ for $i\neq q-1, q$,
  \item ${\cal{H}}^q (K)$ is locally free of finite rank, and
  \item the projection of $K^{q-1}$ onto ${\cal{B}}^q = \im (d)$ locally admits a section.
\end{enumerate}
Then we denote by $\scindage_{[q-1, q]} (K)$ the gerbe of splittings of the complex
\[
  \cdots \to 0 \to K^{q-1} \to K^q \to 0 \to \cdots
\]
concentrated in degrees $0$ and $1$ and with $d$ being equal to the original differential of $K$, rather than $(-1)^{q-1}$ times that; this convention has the consequence that
\[
  \class \scindage_{[q-1, q]} (K) = (-1)^{q-1} \class \scindage (K[q-1])
\]
in $H^2 (X, \underline{\Hom}( {\cal{H}}^q, {\cal{H}}^{q-1}) )$. Combined with Lemma~\ref{lemma:di-scindage-class} and \eqref{eqn:delta-shift-sign}, this implies the following generalization of Lemma~\ref{lemma:di-scindage-class}:
\[ 
    \class \scindage_{[q-1, q]} (K) =  -\delta^q_K.
\]
When there is no confusion as to what $q$ is, we simply write $\scindage (K)$ for $\scindage_{[q-1,q]} (K)$.

\subsection{Truncated Koszul complexes}

Let $K=[K^0\xrightarrow{d} K^1]$ be a two-term complex of modules over $(\mathcal{X}, \mathcal{O})$, and let $q\geq 1$. Using the Koszul complex, one can build another two-term complex, concentrated in degrees $[q-1, q]$:
\[ 
  \tau_{\geq q-1}(\Kos^q(d)) = \left[ \cdots\to 0\to \frac{\exterior^{q-1} K^1 \otimes K^0}{d(\exterior^{q-2}K^1\otimes \Gamma^2 K^0)} \to \exterior^q K^1 \to 0\to \cdots \right].
\]
By Proposition~\ref{prop:alphamaps} we have morphisms
\begin{equation} \label{eqn:alphamaps-Lq}
  \alpha^q\colon  \exterior^q \h^1(K) \To \h^q(\Kos^q(d)) 
  \quad \text{and} \quad  
  \alpha^{q-1}\colon (\exterior^{q-1}\h^1(K))\otimes \h^0(K) \To \h^{q-1}(\Kos^q(d))
\end{equation}
which are isomorphisms if $K^0$, $K^1$, $\h^0(K)$, $\h^1(K)$ are flat. The following result describes the maps $\delta^q_{\Kos^q(d)}$ and hence the truncation $\tau_{\geq q-1}(\Kos^q(d))$.

\begin{proposition} \label{prop:deltaLq}
  Let $K=[K^0\xrightarrow{d} K^1]$ be a two-term complex and let $q\geq 1$. Suppose that $K^0$, $K^1$, $\h^0(K)$, and $\h^1(K)$ are flat. Then the following diagram is commutative
  \[ 
    \xymatrix{
      \exterior^q \h^1(K) \ar[r]^-{\eta^q} \ar[d]_{\alpha^q} & (\exterior^{q-1}\h^1(K))\otimes \h^1(K) \ar[r]^-{\id\otimes \delta^1_K} & (\exterior^{q-1} \h^1(K))\otimes \h^0(K)[2] \ar[d]^{\alpha^{q-1}} \\
      \h^q(\Kos^q(d)) \ar[rr]_{\delta^q_{\Kos^q(d)}} & & \h^{q-1}(\Kos^q(d))[2].
    }
  \]
\end{proposition}

In other words, using the identifications \eqref{eqn:alphamaps-Lq}, we have the equality
\[ 
    \delta^q_{\Kos^q(d)} = (\id\otimes \delta^1_K)\circ \eta^q
\]
of maps $\exterior^q \h^1(K)\to (\exterior^{q-1} \h^1(K))\otimes \h^0(K)[2]$.

\begin{proof}
Let us abbreviate $\h^i(K)$ to $\h^i$. We first check that the two-term complexes $\tau_{\geq q-1}\Kos^q(d)$ and $((\exterior^{q-1} \h^1)\otimes K)[1-q]$ form the middle square inside a commutative diagram with exact rows
\[ 
  \xymatrix@C=1.7em{
    0\ar[r] & (\exterior^{q-1} \h^1) \otimes \h^0 \ar[d]_{\id} \ar[r] & {\displaystyle  \frac{(\exterior^{q-1} K^1)\otimes K^0}{d((\exterior^{q-2} K^1)\otimes \Gamma^2(K^0))} } \ar[d]^\beta \ar[r] & \exterior^q K^1 \ar[d]^\alpha \ar[r] & \exterior^q \h^1 \ar[d]^{\eta^q} \ar[r] & 0\\
    0 \ar[r] & (\exterior^{q-1} \h^1) \otimes \h^0 \ar[r] & (\exterior^{q-1} \h^1)\otimes K^0 \ar[r]_-{\id\otimes (-1)^{q-1}d} & (\exterior^{q-1} \h^1)\otimes K^1 \ar[r] & (\exterior^{q-1} \h^1) \otimes \h^1 \ar[r] & 0.
  }
\]
We define the maps $\alpha$ and $\beta$ as follows. The map $\beta$ is uniquely determined by
\[
  \beta(z_1\wedge\ldots \wedge z_{q-1}\otimes w \,\,\text{ mod } \,\, d(\exterior^{q-2} K^1\otimes \Gamma^2(K^0))) = [z_1]\wedge \ldots\wedge [z_q]\otimes w.
\]
It is well-defined because an elements of the form
\begin{align*} 
  d(z_1\wedge\ldots\wedge z_{q-1}\otimes w^{[2]}) &= z_1\wedge\ldots\wedge z_{q-1}\wedge dw \otimes w, \qquad \text{or} \\
  d(z_1\wedge\ldots\wedge z_{q-1}\otimes wv) &= z_1\wedge\ldots\wedge z_{q-1}\wedge dw \otimes v + z_1\wedge\ldots\wedge z_{q-1}\wedge dv \otimes w.
\end{align*} 
are sent to zero, since $[dw]=0=[dv]$. The map $\alpha$ is the composition
\[
    \exterior^q K^1\xrightarrow{\eta^q} (\exterior^{q-1} K^1)\otimes K^1 \xrightarrow{\text{proj.}\otimes\id} (\exterior^{q-1} \h^1)\otimes K^1.
\]

The commutativity of the left and rightmost squares is trivial to check. To see that the middle square commutes, we take (the class of) $z_1\wedge \ldots\wedge z_{q-1}\otimes w\in (\exterior^{q-1} K^1)\otimes K^0$, and compute
\begin{align*}
  \alpha(d(z_1\wedge \ldots\wedge z_{q-1}\otimes w)) 
  &= \alpha(z_1\wedge \ldots\wedge z_{q-1}\wedge dw) \\
  &= \sum_{i=1}^q (-1)^{i-1} [z_1] \wedge\ldots \widehat{[z_i]}\ldots \wedge [z_{q-1}] \wedge [dw] \otimes z_i + (-1)^{q-1} [z_1]\wedge \ldots \wedge [z_{q-1}] \otimes dw \\
  &=(-1)^{q-1} [z_1]\wedge \ldots \wedge [z_{q-1}]\otimes dw = (-1)^{q-1}(\id\otimes d)(\beta(z_1\wedge \ldots\wedge z_{q-1}\otimes w)).
\end{align*}

Now, thanks to \eqref{eqn:delta-shift-sign}, our commutative diagram of complexes translates into a commutative square in the derived category
\[ 
    \xymatrix{
        \exterior^q \h^1 \ar[r]^-{\delta^q_{\Kos^q(d)}} \ar[d]_{\eta^q} & (\exterior^{q-1} \h^1)\otimes \h^0[2] \ar[d]^\id \\
        (\exterior^q \h^1)\otimes \h^1 \ar[r]_-{\id\otimes\delta^1_K} & (\exterior^{q-1} \h^1)\otimes \h^0[2]
    }
\]
This implies the required assertion.
\end{proof}

\subsection{Two-term truncations of abstract Koszul complexes}

The following result relates the maps $\delta^q_K$ and $\delta^1_K$ for a cdga $K$.
  
\begin{proposition} \label{prop:cdga-deltaq}
  Suppose $2$ is a nonzerodivisor in $\cO$. Let $K$ be a coconnective commutative differential graded algebra such that $K^0$, $\cal{Z}^1 K$, $\h^0(K)$, $\h^1(K)$. Let $q\geq 1$ be an integer such that $\h^{q-1}(K)$ is flat. Then, the following diagram commutes
  \[ 
    \xymatrix{ 
      \exterior^q \h^1(K) \ar[d]_{\text{mult.}} \ar[r]^-{\eta^q} & (\exterior^{q-1} \h^1(K)) \otimes \h^1(K) \ar[r]^-{\id\otimes \delta^1_K} & (\exterior^{q-1} \h^1(K))\otimes \h^0(K) [2]\ar[d]^{\text{mult.}} \\
      \h^q(K) \ar[rr]_{\delta^q_K} & & \h^{q-1}(K)[2].
      }
  \]
\end{proposition}

\begin{proof}
Write $\tau_{\leq 1} K = [K^0\xrightarrow{\partial} \cal{Z}^1 K]$. The proof of Theorem~\ref{thm:koszul-reconstr} (with $m=1$) provides a morphism of complexes 
\[
  \mu\colon \tau_{\geq q-1} \Kos^q(\partial)\To \tau_{[q-1,q]} K.
\]
By functoriality of the maps $\delta^q$, we have a commutative square
\[ 
  \xymatrix@C=1.5em{
    \h^q(\tau_{\geq q-1} \Kos^q(\partial)) \ar[r] \ar[d]_\mu & \h^{q-1}(\tau_{\geq q-1} \Kos^q(\partial))[2] \ar[d] \ar@{}[drr]|-{\displaystyle \text{i.e.}} & & \exterior^q \h^1(K) \ar[r] \ar[d] & (\exterior^{q-1}\h^1(K))\otimes \h^0(K)[2] \ar[d] \\
    \h^q(\tau_{[q-1,q]} K) \ar[r] & \h^{q-1}(\tau_{[q-1,q]} K)[2] & &  \h^q(K) \ar[r] & \h^{q-1}(K)[2].
  }
\]
The assertion then follows from Proposition~\ref{prop:deltaLq}.
\end{proof}

\begin{remark} \label{rmk:achinger-ogus-1}
The proof of \cite[Theorem~4.2.2(1)]{AchingerOgus} implies the assertion of Proposition~\ref{prop:cdga-deltaq} under the stronger assumption that $q!$ is invertible in $\cO$. However, the argument does not use the cdga structure of $K$, only a weaker structure of a commutative monoid in the derived category $K\otimes^\mathbf{L} K\to K$. In particular, the assertion holds for some $E_\infty$-algebras which are not a priori equivalent to cdgas.
\end{remark}

\begin{remark} \label{rmk:achinger-ogus-2}
In \cite[Lemma~2.1.1]{AchingerOgus}, it is shown that the maps $\delta^q_K$ are compatible with the derived tensor product in the following way. If $K$ and $L$ are complexes and $i, j$ are integers such that $\h^i(K)$ and $\h^j(L)$ are flat $\cO$-modules, then the following square commutes.
\[
  \xymatrix{ 
    \h^i(K)\otimes \h^j(L) \ar[d] \ar[rr]^-{\delta^i_K\otimes 1 + (-1)^i\otimes \delta^j_L} & & \h^{i-1}(K)[2]\otimes \h^j(L) \,\,\oplus\,\, \h^i(K)\otimes \h^{j-1}(L)[2] \ar[d] \\
    \h^{i+j}(K\otimes^{\mathbf{L}} L) \ar[rr]_{\delta_{K\otimes^\mathbf{L} L}^{i+j}} & & \h^{i+j-1}(K\otimes^\mathbf{L} L)[2]
  }
\] 
If $q!$ is invertible in $\cO$, so that $\exterior^q \h^1(K)$ is a direct summand of $\h^1(K)^{\otimes q}$, the assertion of Proposition~\ref{prop:cdga-deltaq} can be deduced from this result. 

For illustration, let us see how to do this for $q=2$. We set $L=K$ and $i=j=1$ in the above diagram, obtaining the middle square of the diagram below. 
\[ 
  \xymatrix@R=3em{
    \exterior^2 \h^1(K) \ar[d]_{x\wedge y\,\mapsto\, \frac{1}{2}(x\otimes y - y\otimes x)} \ar[r]^-{\eta^2} & \h^1(K)\otimes\h^1(K) \ar[r]^-{\id\otimes \delta^1_K} & \h^1(K)\otimes \h^0(K)[2] \ar[d]^{\frac 1 2(\mathrm{shuffle}, \id)} \\ 
    \h^1(K)\otimes \h^1(K) \ar[d] \ar[rr]^{\delta^1_K\otimes 1 - 1\otimes \delta^1_K} & & \h^0(K)[2]\otimes \h^1(K) \,\,\oplus\,\, \h^1(K)\otimes \h^0(K)[2] \ar[d] \\
    \h^2(K\otimes^\mathbf{L} K) \ar[d] \ar[rr]_{\delta^2_{K\otimes^\mathbf{L} K}} & & \h^1(K\otimes^\mathbf{L} K)[2] \ar[d] \\
    \h^2(K) \ar[rr]_{\delta^2_K} & & \h^1(K)[2]
  }
\]
Here, the bottom square certifies the functoriality of $\delta^2$ with respect to the multiplication map $K\otimes^\mathbf{L} K\to K$. Commutativity of the top square is easy to check. Then, commutativity of the exterior square gives the required assertion.
\end{remark}

\begin{corollary} \label{cor:two-step1}
  Suppose $2$ is a nonzerodivisor in $\cO$. Let $K$ be an abstract Koszul complex satisfying the flatness condition \eqref{eqn:flatness} and let $q\geq 1$. We have the following commutative diagram
  \[ 
    \xymatrix{ 
      \exterior^q \h^1(K) \ar[d]_\vsim \ar[r]^-{\eta^q} & (\exterior^{q-1} \h^1(K)) \otimes \h^1(K) \ar[r]^-{\id\otimes\delta^1_K} & \exterior^{q-1} \h^1(K)[2]\ar[d]^\vsim \\
      \h^q(K) \ar[rr]_{\delta^q_K} & & \h^{q-1}(K)[2].
      }
  \]
  In other words, using the vertical identifications, we have the equality
  \[ 
    \delta^q_K = (\id\otimes \delta^1_K)\circ \eta^q
  \]
  in $\Hom(\h^q(K), \h^{q-1}(K)[2])$.
\end{corollary}

\begin{corollary} \label{cor:clsc-formula}
  Suppose $2$ is a nonzerodivisor in $\cO$. Let $K$ be an abstract Koszul complex satisfying the flatness condition \eqref{eqn:flatness} and let $q\geq 1$. Then
  \[ 
    \class \scindage_{[q-1,q]}(K) = \eta^q(\class \scindage (\tau_{\leq 1} K)) .
  \]
\end{corollary}

\begin{corollary} \label{cor:d2-formula}
  Suppose that $2$ is a nonzerodivisor in $\cO$. Let $K$ be an abstract Koszul complex satisfying the flatness condition \eqref{eqn:flatness}. Then the differential 
  \[ 
    d_2^{pq} \colon H^p(X, \h^q(K)) \to H^{p+2}(X, \h^{q-1}(K))
  \]
  equals the cup product with the class 
  \[
    \delta^1_K = -\class \scindage(\tau_{\leq 1}K) \in H^2(X, \underline{\Hom}(\h^1(K), \h^0(K)))
  \]
  followed by evaluation.
\end{corollary}

\subsection{Morphisms of gerbes of splittings}
\label{ss:construction}

Let $K$ be an abstract Koszul complex satisfying the flatness condition \eqref{eqn:flatness}. In Corollary~\ref{cor:clsc-formula}, under the assumption that $2$ is a nonzerodivisor in $\cO$, we calculated the gerbe classes $\class \scindage_{[q-1,q]} K$ in terms of the class $\class\scindage (\tau_{\leq 1} K)$. Below, we promote this equality into a morphism of gerbes.

\begin{theorem} \label{thm-s2.2.1}
 For each integer $q\ge 1$, there is a morphism
 $$
 \wedge^q : \scindage ( \tau_{\leq 1} K ) \To \scindage ( \tau_{[q-1, q]} K )
 $$
 of gerbes over $X$, under which the obstruction classes correspond by the relation
 $$
 \class \scindage (\tau_{[q-1, q]} K ) = \ctr^q ( \class \scindage \tau_{\leq 1} (K)),
 $$
 where $\ctr^q: \underline{\Hom}( {\cal{H}}^1, {\cal{H}}^0 ) \to \underline{\Hom}( {\cal{H}}^q, {\cal{H}}^{q-1} )$ denotes the morphism which maps a local section $f$ of the source to the one of the target by the formula
 $$
 \ctr^q (f): \omega_1 \wedge \cdots \wedge \omega_q \mapsto
    \sum_{j=1}^{q} (-1)^{j-1} f(\omega_j) \omega_1 \wedge \cdots \wedge \omega_{j-1} \wedge \omega_{j+1} \wedge \cdots \wedge \omega_q.
 $$
\end{theorem}
\noindent (Compare the formula for $\ctr^q$ with the explicit formula for $\eta^q$ in Section \ref{section-Koszul-complexes}.)

\medskip

\noindent {\bf Notations.} Before proceeding to the proof, we gather some notations concerning \v{C}ech cohomology. We denote by $\check{\cal{C}}(U_{\bullet}, K^{\bullet})$
the {\v{C}}ech resolution of a complex $K^{\bullet}$ with respect to a hypercovering $U_{\bullet}$.
The differential induced by that of $K^{\bullet}$ will still be denoted by $d$,
while the {\v{C}}ech differential on the component
$\check{\cal{C}}(U_p, K^q)$: 
$$
(-1)^q \sum_{i=0}^{p+1} (-1)^i d_i^{\ast}
$$
will be denoted by $\check{d}$. Then the total differential
$$
D = d + \check{d}
$$
is the differential of the total complex $\check{\cal{C}}(U_{\bullet}, K^{\bullet})$.

When we compute the obstruction classes, we will use some notations which may not be standard. 
As usual, for each integer $m \ge -1$, we denote by $[m]$ the set of integers $i$ such that $0 \le i \le m$ (empty set for $[-1]$). 
And we denote by 
$d_{ij} : [m-2] \to [m]$
the unique increasing injection omitting $i$ and $j$, where $0\le i<j\le m$.
For example, for $m = 2$, we have
$$
d_{02} = d_2 \circ d_0 = d_0 \circ d_1 : [0] \to [2]
$$
(which maps $0$ onto $1$), where $d_i : [m-1] \to [m]$ denotes the unique increasing injection omitting $i$.

On the other hand, we denote by $\prj_i : [0] \to [m]$ (resp. $\prj_{ij} : [1] \to [m]$) 
the unique map sending $0$ to $i$ (resp. $0$ to $i$ and $1$ to $j$), for $0\le i \le m$ (resp. for $0\le i<j\le m$).

\begin{proof}[Proof of Theorem~\ref{thm-s2.2.1}]

In order to prove Theorem~\ref{thm-s2.2.1}, we first describe the morphism, show that it is well-defined, and then calculate the obstruction class.

\medskip

\noindent {\bf Construction of the functor $\wedge^q$.} We construct $\wedge^q : \scindage \tau_{\leq 1} K \to \scindage \tau_{[q-1, q]} K$ by stackifying a morphism between the corresponding prestacks: 
$\scindage' \tau_{\leq 1} K \to \scindage' \tau_{[q-1, q]} K$.

Given an object of $\scindage' (\tau_{\leq 1} K)$ over $U$, that is, a section $s : {\cal{H}}^1 \to {\cal{Z}}^1$ of the projection ${\cal{Z}}^1 \to {\cal{H}}^1$ over $U$, we define $\wedge^q (s)$ as the composite morphism
$$
\xymatrix{
{\cal{H}}^q \simeq \wedge^q ( {\cal{H}}^1 ) \ar[r]^-{\wedge^q (s)} & \wedge^q ( {\cal{Z}}^1 ) \ar[r]^-{\product} & {\cal{Z}}^q,
}
$$
where $\product$ means product; it is clearly a section of ${\cal{Z}}^q \to {\cal{H}}^q$ over $U$.

Let $s_0$ and $s_1$ be two objects of $\scindage' (\tau_{\leq 1} K)$ over $U$ and let $h$ be a homotopy from $s_0$ to $s_1$. Then we need to define a corresponding homotopy $\wedge^q (h)$ from $\wedge^q (s_0)$ to $\wedge^q (s_1)$. We first define a map
$$
\wedge'^q (h) : ( {\cal{H}}^1 )^{\otimes q} \To K^{q-1} / {\cal{B}}^{q-1}
$$
by letting it send $\omega_1 \otimes \cdots \otimes \omega_q$ (where $\omega_j$ are local sections of ${\cal{H}}^1$) to the class of
$$
\sum_{j=1}^{q} (-1)^{j-1} h(\omega_j ) s_0 (\omega_1) \wedge \cdots \wedge s_0 (\omega_{j-1} ) \wedge s_1 (\omega_{j+1}) \wedge \cdots \wedge s_1 (\omega_q)
$$
modulo ${\cal{B}}^{q-1}$.

It is easy to show that it factors through $\wedge^q {\cal{H}}^1$ : If, say, $\omega_1 = \omega_2 = \omega$, then the alternating sum on the right reduces to the difference of the first two terms
$$
h(\omega) s_1(\omega) \wedge s_1 (\omega_3) \wedge \cdots \wedge s_1 (\omega_q) 
- h(\omega) s_0 (\omega) \wedge s_1 (\omega_3) \wedge \cdots \wedge s_1 (\omega_q) ,
$$
which is equal to
$$
h(\omega) dh(\omega) \wedge s_1 (\omega_3) \wedge \cdots \wedge s_1 (\omega_q)
$$
which is a coboundary {\bf since $2$ is invertible}. Thus we defined $\wedge^q (h)$.
$$
\xymatrix{
({\cal{H}}^1)^{\otimes q} \ar[rd]^-{\wedge'^q (h)} \ar@{>>}[d] \\
\wedge^q ( {\cal{H}}^1 ) \ar@{-->}[r]_-{\wedge^q (h) } & K^{q-1} / {\cal{B}}^{q-1}
}
$$
Then the following calculation shows that $\wedge^q (h)$ is really a homotopy:
\begin{align*}
  \left( \wedge^q (s_1)  -  \wedge^q (s_0) \right) ( \omega_1 \wedge \cdots \wedge \omega_q ) 
  & = \sum_{j=1}^q s_0 ( \omega_1 ) \wedge \cdots \wedge s_0 ( \omega_{j-1} ) 
    \wedge \{ s_1 (\omega_j ) - s_0 (\omega_j ) \} \wedge s_1 (\omega_{j+1} ) \wedge \cdots \wedge s_1 (\omega_q ) \\
  & = \sum_{j=1}^q s_0 ( \omega_1 ) \wedge \cdots \wedge s_0 ( \omega_{j-1} ) 
    \wedge \{ dh(\omega_j ) \} \wedge s_1 (\omega_{j+1} ) \wedge \cdots \wedge s_1 (\omega_q ) \\
  & = \sum_{j=1}^q (-1)^{j-1} dh(\omega_j ) \wedge s_0 ( \omega_1 ) \wedge \cdots \wedge s_0 ( \omega_{j-1} ) \wedge s_1 (\omega_{j+1} ) \wedge \cdots \wedge s_1 (\omega_q ) \\
  & =  d [ \wedge^q (h) (\omega_1 \wedge \cdots \wedge \omega_q ) ].
\end{align*}

\noindent {\bf Functoriality of $\wedge^q$.} 
Now in order to show that the morphism $\wedge^q$ is a functor, we must show that it is compatible with the composition of homotopies; so let $h: s_0 \Rightarrow s_1$ and $h': s_1 \Rightarrow s_2$ be two such in the source. We first define a second homotopy operator:
\begin{eqnarray*}
 H_2^q (h, h') : ( {\cal{H}}^1 )^{\otimes q} & \to & K^{q-2} \\
 \omega_1 \otimes \cdots \otimes \omega_q & \mapsto & \sum_{1\le j < k \le q} (-1)^{j+k+1} h(\omega_j ) h'(\omega_k ) s(j,k),
\end{eqnarray*}
where $s(j, k)$ is equal to
$$
s_0 (\omega_1 ) \wedge \cdots \wedge s_0 (\omega_{j-1} ) \wedge
    s_1 (\omega_{j+1} ) \wedge \cdots \wedge s_1 (\omega_{k-1} ) \wedge
        s_2 (\omega_{k+1} ) \wedge \cdots \wedge s_2 (\omega_q ).
$$
To show that $\wedge^q (h + h')$ and $\wedge^q(h) + \wedge^q(h')$ are the same homotopies, it suffices to
demonstrate the formula
$$
[ \wedge^q (h+h') - \{ \wedge^q (h) + \wedge^q (h') \} ] (\omega_1 \wedge \cdots \wedge \omega_q )
  = d H_2^q (\omega_1 \otimes \cdots \otimes \omega_q ).
$$
One expands the left hand side and groups the terms involving $h$ and $h'$ separately:
\begin{eqnarray*}
 & & \sum_{j=1}^q (-1)^{j-1} \{ h(\omega_j ) + h'(\omega_j ) \} s_0 (\omega_1) \wedge \cdots \wedge s_0 (\omega_{j-1} ) \wedge s_2 (\omega_{j+1} ) \wedge \cdots \wedge s_2 (\omega_q ) \\
 & & \,\, - \sum_{j=1}^q (-1)^{j-1} h(\omega_j ) s_0 (\omega_1) \wedge \cdots \wedge s_0 (\omega_{j-1} ) \wedge s_1 (\omega_{j+1} ) \wedge \cdots \wedge s_1 (\omega_q ) \\
 & & \,\, - \sum_{j=1}^q (-1)^{j-1} h'(\omega_j ) s_1 (\omega_1) \wedge \cdots \wedge s_1 (\omega_{j-1} ) \wedge s_2 (\omega_{j+1} ) \wedge \cdots \wedge s_2 (\omega_q ) \\
 & = & \sum_{j=1}^q (-1)^{j-1} h(\omega_j ) s_0 (\omega_1 ) \wedge \cdots \wedge s_0 (\omega_{j-1} )
 \wedge \{ s_2 (\omega_{j+1} ) \wedge \cdots \wedge s_2 (\omega_q ) - s_1 (\omega_{j+1} ) \wedge \cdots \wedge s_1 (\omega_q ) \} \\
 & + & \sum_{k=1}^q (-1)^{k-1} h' (\omega_k ) \{ s_0 (\omega_1 ) \wedge \cdots \wedge s_0 (\omega_{k-1} ) - s_1 (\omega_1 ) \wedge \cdots \wedge s_1 (\omega_{k-1} ) \} \wedge s_2 (\omega_{k+1} ) \wedge \cdots \wedge s_2 (\omega_q).
\end{eqnarray*}
Now the differences in the curly brackets are themselves alternating sums, so
\begin{eqnarray*}
 & = & \sum_{j=1}^q (-1)^{j-1} h(\omega_j ) s_0 (\omega_1 ) \wedge \cdots \wedge s_0 (\omega_{j-1} ) \wedge \\
 & & \,\,\,\, \wedge \left\{ \sum_{k>j} (-1)^{k-(j+1)} dh' (\omega_k ) \wedge s_1 (\omega_{j+1}) \wedge \cdots \wedge s_1 (\omega_{k-1} ) \wedge s_2(\omega_{k+1} ) \wedge \cdots \wedge s_2 (\omega_q ) \right\} \\
 & + & \sum_{k=1}^q (-1)^{k-1} h' (\omega_k ) \left\{ \sum_{j<k} (-1)^{j-1} (-d h(\omega_j )) \wedge
    s_0 (\omega_1 ) \wedge \cdots \wedge s_0 ( \omega_{j-1} ) \wedge s_1 (\omega_{j+1} ) \wedge \cdots s_1 (\omega_{k-1}) \right\} \wedge \\
 & & \,\,\,\, \wedge s_2 (\omega_{k+1} ) \wedge \cdots \wedge s_2 (\omega_q) \\
 & = & \sum_{1\le j < k \le q} (-1)^{j+k+1} \{ h(\omega_j ) dh' (\omega_k ) + h'(\omega_k ) dh(\omega_j )\} \wedge s(j,k),
\end{eqnarray*}
and this is now equal to:
\begin{equation}
\label{eqn-second-htpy}
  = d H_2^q (h, h') (\omega_1 \otimes \cdots \otimes \omega_q ).
\end{equation}
This completes the proof of the fact that $\wedge^q$ is a functor.

\medskip

\noindent {\bf Calculation of obstruction classes.}
Finally, we relate the obstruction elements. Let $U_{\bullet} \to X$ be an open hypercovering such that one has
\begin{enumerate}
 \item A section $s : {\cal{H}}^1 \to {\cal{Z}}^1$
 of the canonical projection ${\cal{Z}}^1 \to {\cal{H}}^1$ over $U_0$ and
 \item A homotopy $h : {\cal{H}}^1 \to K^0$ over $U_1$ satisfying
 \begin{equation}
  \label{eqn-s2.1}
  d_1^{\ast} s - d_0^{\ast} s = d h
 \end{equation}
\end{enumerate}
Then by definition, the class of
$$
\obs = \obs_1 = d_0^{\ast} h - d_1^{\ast} h + d_2^{\ast} h \in \Gamma ( U_2, \underline{\Hom} ( {\cal{H}}^1, {\cal{H}}^0 ) )
$$
in $H^2 (X, \underline{\Hom}( {\cal{H}}^1, {\cal{H}}^0 ))$ is $\class \scindage \tau_{\le 1} K$.
On the other hand, by applying $\wedge^q$ to $s$ and $h$, one sees that the class of
$$
\obs_q = d_0^{\ast} (\wedge^q h ) - d_1^{\ast} (\wedge^q h )  + d_2^{\ast} (\wedge^q h)
\in \Gamma ( U_2, \underline{\Hom} ( {\cal{H}}^q, {\cal{H}}^{q-1} ) )
$$
in $H^2 (X, \underline{\Hom}( {\cal{H}}^q, {\cal{H}}^{q-1} ))$ is $\class \scindage \tau_{[q-1, q]} K$.

Now let $\omega_1 \wedge \cdots \wedge \omega_q$ be a local section of ${\cal{H}}^q \simeq \wedge^q {\cal{H}}^1$. 
Then the evaluation of $\obs_q$ at $\omega_1 \wedge \cdots \wedge \omega_q$ is equal to
\begin{eqnarray*}
 & & \sum_{j=1}^q (-1)^{j+1} d_0^{\ast} h(\omega_j ) d_{01}^{\ast} s(\omega_1 ) \wedge \cdots \wedge d_{01}^{\ast} s(\omega_{j-1} ) \wedge d_{02}^{\ast} s (\omega_{j+1} ) \wedge \cdots d_{02}^{\ast} s(\omega_q ) \\
 & - & \sum_{j=1}^q (-1)^{j+1} d_1^{\ast} h(\omega_j ) d_{01}^{\ast} s(\omega_1 ) \wedge \cdots \wedge d_{01}^{\ast} s(\omega_{j-1} ) \wedge d_{12}^{\ast} s (\omega_{j+1} ) \wedge \cdots d_{12}^{\ast} s(\omega_q ) \\
 & + & \sum_{j=1}^q (-1)^{j+1} d_2^{\ast} h(\omega_j ) d_{02}^{\ast} s(\omega_1 ) \wedge \cdots \wedge d_{02}^{\ast} s(\omega_{j-1} ) \wedge d_{12}^{\ast} s (\omega_{j+1} ) \wedge \cdots d_{12}^{\ast} s(\omega_q ).
\end{eqnarray*}
One groups the terms around the second sum and gets
\begin{eqnarray*}
 & & \sum_{j=1}^q (-1)^{j+1} (d_0^{\ast} h - d_1^{\ast} h + d_2^{\ast} h ) (\omega_j ) \cdot \\
 & & \,\, \cdot d_{01}^{\ast} s(\omega_1 ) \wedge \cdots \wedge d_{01}^{\ast} s(\omega_{j-1} ) \wedge d_{12}^{\ast} s( \omega_{j+1} ) \wedge \cdots \wedge d_{12}^{\ast} s(\omega_q ) \\
 & + & \sum_{j=1}^q (-1)^{j+1} d_0^{\ast} h (\omega_j ) d_{01}^{\ast} s(\omega_1 ) \wedge \cdots \wedge d_{01}^{\ast} s(\omega_{j-1} ) \wedge \\
 & & \,\, \wedge \{ d_{02}^{\ast} s (\omega_{j+1} ) \wedge \cdots \wedge d_{02}^{\ast} s(\omega_q ) - d_{12}^{\ast} s (\omega_{j+1} ) \wedge \cdots \wedge d_{12}^{\ast} s(\omega_q ) \} \\
 & + & \sum_{j=1}^q (-1)^{j+1} d_2^{\ast} h(\omega_j ) \{ d_{02}^{\ast} s( \omega_1 ) \wedge \cdots \wedge d_{02}^{\ast} s(\omega_{j-1} ) - d_{01}^{\ast} s( \omega_1 ) \wedge \cdots \wedge d_{01}^{\ast} s(\omega_{j-1} ) \} \wedge \\
 & & \,\, \wedge d_{12}^{\ast} s(\omega_{j+1} ) \wedge \cdots \wedge d_{12}^{\ast} s(\omega_q ).
\end{eqnarray*}
The first alternating sum reduces to the ``main'' term we want, when taken modulo the coboundaries (${\cal{B}}^{q-1}$). In the last two sums, one first notes that, as $h$ is a homotopy from $d_0^{\ast} s$ to $d_1^{\ast} s$, it follows that $d_0^{\ast} h$ is a homotopy from $d_0^{\ast} d_0^{\ast} s$ to $d_0^{\ast} d_1^{\ast} s$, that is,
$$
d_0^{\ast} h : d_{01}^{\ast} s \Rightarrow d_{02}^{\ast} s.
$$
Similarly, $d_2^{\ast} h$ is a homotopy from $d_{02}^{\ast} s$ to $d_{12}^{\ast} s$.
Essentially by repeating the last three equalities leading up to (\ref{eqn-second-htpy}),
this time with a minus sign, one sees that the last two sums add up to
$$
-d H^q_2 (d_0^{\ast} h, d_2^{\ast} h ) (\omega_1 \otimes \cdots \omega_q ),
$$
which is a coboundary. Therefore, reducing modulo ${\cal{B}}^{q-1}$, one gets
$$
\eval (\obs_q , \omega_1 \wedge \cdots \wedge \omega_q ) = \eval (\ctr^q (\obs_1 ),  \omega_1 \wedge \cdots \wedge \omega_q ).
$$
This means
$$
\class \scindage \tau_{[q-1, q]} K = \ctr^q \class \scindage \tau_{\leq 1} K,
$$
which completes the proof.
\end{proof}

\begin{remark}
The construction of the map between gerbes can also be carried out using the language of higher topos theory \cite{Lurie}. Let us give a brief outline.

Let $p\colon Y\to Z$ be a map in the homotopy category of spaces, or more generally in any $\infty$-category $\cal{C}$. On can then build the space $\scindage(p)$ of splittings of $p$ as the homotopy fiber of 
\[
    p\colon \Hom_{\cal{C}}(Z, Y)\To \Hom_{\cal{C}}(Z, Z).
\]
over the identity ${\rm id}_Z$. Similarly, if $p\colon Y\to Z$ is a map in the derived $\infty$-category of a ringed topos $(X, \cO)$, one obtains a sheaf of spaces $\underline{\scindage}(p)$ of splittings of $p$. 

In the special case when $Y$ is a two-term complex $K=[K^0\xrightarrow{d} K^1]$ satisfying the conditions in \S\ref{ss:gerbe-of-splittings} and $Y\to Z$ is the projection $K\to \h^1(K)[-1]$, then $\underline{\scindage}(p)$ is a sheaf of groupoids (a stack) and can be identified with the gerbe of splittings $\scindage K$. 

Applying the functor $\tau_{\geq q-1} L\Gamma^q$ to the map $p\colon K\to \h^1(K)[-1]$ one obtains (simply by functoriality) a~morphism of sheaves of spaces
\[ 
    \underline{\scindage}(p)\To \underline{\scindage}(\tau_{\geq q-1} L\Gamma^q(p)). 
\]
By inspection, the map $\tau_{\geq q-1} L\Gamma^q(p)$ is the projection 
\[
    \tau_{[q-1,q]}\Kos^q(d)\To \h^q(\Kos^q(d))[-q] = \exterior^q\h^1(K)[-q].
\]
This way one obtains by abstract nonsense a morphism of gerbes $\scindage K\to \scindage(\tau_{[q-1,q]}\Kos^q(d))$.

In the case when $K$ is an abstract Koszul complex satisfying the flatness condition \eqref{eqn:flatness}, the morphism of gerbes $\scindage \tau_{\leq 1} K\to \scindage \tau_{[q-1,q]} K$ obtained this way should agree with the one constructed in Theorem~\ref{thm-s2.2.1}, though we did not check it.
\end{remark}

\section{Gerbes of splittings of the de Rham complex}
\label{s:gerbe-calculation}

Our method of explicating the truncations $\tau_{[q-1, q]} K$ for an abstract Koszul complex $K$ in terms of the truncation $\tau_{\leq 1} K$ requires that $2$ is a nonzerodivisor. In this section, we describe these two-term truncations in the case of the de Rham complex in characteristic $p>0$ by calculating the class
$$
\class \scindage \left( \tau_{[q-1, q]} F_{\ast} \Omega_{X/S}^{\bullet} \right)
$$
The calculation uses more information about the de Rham complex than it being an abstract Koszul complex, namely the nature of the Cartier isomorphism (which we use only for $p=2$). As a corollary, we deduce that $\tau_{[q-1, q]} (F_*\Omega^\bullet_{X/S})$ is decomposable if $\tau_{\leq 1}(F_* \Omega^\bullet_{X/S})$ is, and obtain a description of the $d_2$ differentials in the conjugate spectral sequence. 

\begin{theorem} \label{thm-s2.3.1}
 Let $S$ be a scheme of characteristic $p > 0$ and $X/S$ a smooth separated scheme of finite type. Then for each integer $q$, the class
 $$
 \class \scindage ( \tau_{[q-1, q]} F_{\ast} \Omega^{\bullet}_{X/S} )
 $$
 is the image of the class
 $$
 \class \scindage ( \tau_{\le 1} F_{\ast} \Omega^{\bullet}_{X/S} )
 $$
 under the contraction map (described in Theorem \ref{thm-s2.2.1}).
\end{theorem}

\begin{proof}
We put $K = F_* \Omega_{X/S}^{\bullet}$, with $F : X \to X'$ the relative Frobenius of $X/S$.

To calculate the class, we take an open hypercovering $U_{\bullet}\to X'$ such that
\begin{enumerate}
 \item Over $U_0$, one has a section $s : {\cal{H}}^1 \to {\cal{Z}}^1$ of the projection ${\cal{Z}}^1 \to {\cal{H}}^1$ and a section
 $$
 \sigma^{(q)} : {\cal{H}}^q \to ( {\cal{H}}^1)^{\otimes q}
 $$
 of the canonical projection $({\cal{H}}^1)^{\otimes q} \to {\cal{H}}^q$; and
 \item Over $U_1$, one has a homotopy $h : {\cal{H}}^1 \to K^0$ such that
 $$
 dh = d_1^{\ast} s - d_0^{\ast} s : {\cal{H}}^1 \to  {\cal{Z}}^1.
 $$
\end{enumerate}
(Let us recall that ${\cal{H}}^1$, and hence ${\cal{H}}^q = \wedge^q {\cal{H}}^1$ for all integers $q$, are locally free over ${\cal{O}}_{X'}$.) 
The locally free kernel of the projection $( {\cal{H}}^1 )^{\otimes q} \to {\cal{H}}^q$ being denoted by ${\cal{I}}^q$, the $1$-cocycle
$$
d_0^{\ast} \sigma^{(q)} - d_1^{\ast} \sigma^{(q)} \in \Gamma \left( U_1, \underline{\Hom}_{{\cal{O}}_{X'}} ( {\cal{H}}^q, {\cal{I}}^q ) \right)
$$
represents the obstruction, in $H^1 ( X', \underline{\Hom} ( {\cal{H}}^q , {\cal{I}}^q )) = \Ext^1_{{\cal{O}}_{X'}} ( {\cal{H}}^q , {\cal{I}}^q )$, to the global existence of a section.

Let us calculate the class
$$
\class \scindage \tau_{[q-1, q]} K \in H^2 (X', \underline{\Hom} ( {\cal{H}}^q, {\cal{H}}^{q-1} ) )
$$
in characteristic $p\ge 2$. For ease of notation, we denote $\sigma^{(q)}$ simply by $\sigma$ when no confusion is likely.

To do so, we may choose the composite morphism
$$
\xymatrix{
{\cal{H}}^q \ar[r]^-{\sigma^{(q)}} & ( {\cal{H}}^1 )^{\otimes q} \ar[r]^-{s^{\otimes q}} &
        ( {\cal{Z}}^1 )^{\otimes q} \ar[r]^-{\wedge} & {\cal{Z}}^q
}
$$
which we denote by $(s^{\wedge q}) \circ \sigma$, as the section of the projection ${\cal{Z}}^q \to {\cal{H}}^q$ over $U_0$.

Then one forms (the negative of) the \v{C}ech difference
\begin{eqnarray*}
 & & d_1^{\ast}(s^{\wedge q}) \circ d_1^{\ast} \sigma - d_0^{\ast} (s^{\wedge q}) \circ d_0^{\ast} \sigma \\
 & = & [ (d_1^{\ast} s)^{\wedge q} - (d_0^{\ast} s)^{\wedge q} ] \circ d_0^{\ast} \sigma - (d_1^{\ast} s)^{\wedge q} \circ ( d_0^{\ast} \sigma - d_1^{\ast} \sigma ).
\end{eqnarray*}
One notes that the second term is zero, since the image of $d_0^{\ast} \sigma - d_1^{\ast} \sigma$ is contained in ${\cal{I}}^q$, which in turn is annihilated by $(d_1^{\ast} s)^{\wedge q}$, for the wedge product is strictly graded commutative.

Then one expresses the remaining first term as the differential of something:
\begin{eqnarray*}
& & ((d_1^{\ast} s)^{\wedge q} - (d_0^{\ast} s)^{\wedge q} ) (\omega_1 \otimes \cdots \otimes \omega_q ) \\
& = & \sum_{j=1}^{q} (d_0^{\ast} s) \omega_1 \wedge \cdots \wedge (d_0^{\ast} s) \omega_{j-1} \wedge 
    (d_1^{\ast} s - d_0^{\ast} s) \omega_j \wedge (d_1^{\ast} s) \omega_{j+1} \wedge \cdots \wedge (d_1^{\ast} s) \omega_q \\
& = & d \sum_{j=1}^{q} (-1)^{j+1} h(\omega_j) (d_0^{\ast} s) \omega_1 \wedge \cdots \wedge (d_0^{\ast} s) \omega_{j-1} \wedge (d_1^{\ast} s) \omega_{j+1} \wedge \cdots \wedge (d_1^{\ast} s) \omega_q.
\end{eqnarray*}
One defines $\eta = \eta^{(q)} = \eta (\omega_1 \otimes \cdots \otimes \omega_q )$ to be
$$
\sum_{j=1}^{q} (-1)^{j+1} h(\omega_j) (d_0^{\ast} s) \omega_1 \wedge \cdots \wedge (d_0^{\ast} s) \omega_{j-1} \wedge (d_1^{\ast} s) \omega_{j+1} \wedge \cdots \wedge (d_1^{\ast} s) \omega_q
$$
in order to have a commutative diagram
$$
\xymatrix{
{\cal{H}}^q \ar[rrr]^-{d_1^{\ast} (s^{\wedge q}) \circ d_1^{\ast} \sigma - d_0^{\ast}(s^{\wedge q}) \circ d_0^{\ast} \sigma } \ar[rd]_-{d_0^{\ast} \sigma} & & & {\cal{Z}}^q \\
& ( {\cal{H}}^1 )^{\otimes q} \ar[r]^-{\overline{\eta}} & K^{q-1}/{\cal{B}}^{q-1} \ar[ur]_-{d} & ,
}
$$
in which $\overline{\eta}$ means the composite of $\eta$ followed by $K^{q-1} \to K^{q-1} / {\cal{B}}^{q-1}$.

With this, we calculate the class of the gerbe by forming the \v{C}ech difference
\begin{eqnarray}
(d_0^{\ast} - d_1^{\ast} + d_2^{\ast} ) (\overline{\eta} \circ d_0^{\ast} \sigma) 
& = & d_0^{\ast} \overline{\eta} \circ d_{01}^{\ast} \sigma 
    - d_1^{\ast} \overline{\eta} \circ d_{01}^{\ast} \sigma
    + d_2^{\ast} \overline{\eta} \circ d_{02}^{\ast} \sigma \\
& = & (d_0^{\ast} - d_1^{\ast} + d_2^{\ast} ) \overline{\eta} \circ d_{01}^{\ast} \sigma 
    - d_2^{\ast} \overline{\eta} \circ (d_{01}^{\ast} \sigma - d_{02}^{\ast} \sigma ).
 \label{eqn-s2.4}
\end{eqnarray}

Let us put $\obs_1 = (d_0^{\ast} - d_1^{\ast} + d_2^{\ast}) h$, which represents the class of $\scindage \tau_{\le 1} K$. Then the first summand of (\ref{eqn-s2.4}) can be expressed in terms of $\obs_1$:
\begin{eqnarray*}
 & & (d_0^{\ast} - d_1^{\ast} + d_2^{\ast}) \overline{\eta} (\omega_1 \otimes \cdots \otimes \omega_q ) \\
 & = & \sum_{j=1}^q (-1)^{j+1} d_0^{\ast} h (\omega_j) d_{01}^{\ast} s(\omega_1 ) \wedge \cdots \wedge d_{01}^{\ast} s(\omega_{j-1} )   \wedge   d_{02}^{\ast} s (\omega_{j+1} ) \wedge \cdots \wedge d_{02}^{\ast} s (\omega_q ) \\
 & & - \sum_{j=1}^q (-1)^{j+1} d_1^{\ast} h (\omega_j) d_{01}^{\ast} s(\omega_1 ) \wedge \cdots \wedge d_{01}^{\ast} s(\omega_{j-1})   \wedge   d_{12}^{\ast} s (\omega_{j+1} ) \wedge \cdots \wedge d_{12}^{\ast} s (\omega_q ) \\
 & & + \sum_{j=1}^q (-1)^{j+1} d_2^{\ast} h (\omega_j) d_{02}^{\ast} s(\omega_1 ) \wedge \cdots \wedge d_{02}^{\ast} s(\omega_{j-1})   \wedge   d_{12}^{\ast} s (\omega_{j+1} ) \wedge \cdots \wedge d_{12}^{\ast} s (\omega_q ) \\
 & = & \sum_{j=1}^q (-1)^{j+1} \obs_1 (\omega_1) d_{01}^{\ast} s(\omega_1 ) \wedge \cdots \wedge d_{01}^{\ast} s(\omega_{j-1} ) \wedge d_{12}^{\ast} s(\omega_{j+1} ) \wedge \cdots \wedge d_{12}^{\ast} s(\omega_q ) \\
 & & + \sum_{j=1}^q (-1)^{j+1} d_0^{\ast} h(\omega_j ) d_{01}^{\ast} s(\omega_1) \wedge \cdots \wedge d_{01}^{\ast} s(\omega_{j-1} ) \wedge \left\{ (d_{02}^{\ast} s^{\wedge (q-j)} -  d_{12}^{\ast} s^{\wedge (q-j)}) (\omega_{j+1} \otimes \cdots \otimes \omega_q ) \right\} \\
 & & + \sum_{j=1}^q (-1)^{j+1} d_2^{\ast} h(\omega_j ) \left\{ (d_{02}^{\ast} s^{\wedge (j-1)} - d_{01}^{\ast} s^{\wedge (j-1)} ) (\omega_1 \otimes \cdots \otimes \omega_{j-1}) \right\}
 \wedge d_{12}^{\ast} s(\omega_{j+1} ) \wedge \cdots \wedge s(\omega_q ).
\end{eqnarray*}
Again as in the three equalities leading up to (\ref{eqn-second-htpy}), the differences in the curly brackets are
themselves alternating sums, and one sees that the sum of the last two alternating sums is equal to
$$
- d H_2^q ( d_0^{\ast} h, d_2^{\ast} h ) (\omega_1 \otimes \cdots \otimes \omega_q ),
$$
hence is zero modulo ${\cal{B}}^{q-1}$. On the other hand, the first alternating sum is equal to
$$
\eval ( \ctr (\obs_1) , \omega_1 \wedge \cdots \wedge \omega_q ).
$$

Now we analyze the second summand in (\ref{eqn-s2.4}). It is the cup product of two cohomology classes:
\begin{eqnarray*}
 \overline{\eta}|_{{\cal{I}}^{q}} \in \Gamma (U_1, \underline{\Hom}( {\cal{I}}^q, {\cal{H}}^{q-1} ) )
    & \mbox{ representing } 
    & [\overline{\eta}|_{{\cal{I}}^q} ] \in \Ext^1 ({\cal{I}}^q, {\cal{H}}^{q-1} ) \,\, \mbox{ and } \\
  d_0^{\ast} \sigma - d_1^{\ast} \sigma \in \Gamma(U_1, \underline{\Hom}( {\cal{H}}^q, {\cal{I}}^{q} ) )
    & \mbox{ representing }
    & [\sigma] \in \Ext^1 ( {\cal{H}}^q, {\cal{I}}^{q} ).
\end{eqnarray*}
When $q$ is less than $p = \mathrm{char.} \, (S)$, $[\sigma]$ is zero, for in this case one disposes of a canonical section of 
\[
    ({\cal{H}}^1 )^{\otimes q} \to {\cal{H}}^q,
\]
namely the anti-symmetrization.

On the other hand, if $p$ is odd, then $[\overline{\eta}|_{{\cal{I}}^q} ]$ is zero, because (even more strongly) $\overline{\eta}$ itself kills ${\cal{I}}^q$: for example, it maps a local section
$$
\omega \otimes \omega \otimes \omega_3 \otimes \cdots \otimes \omega_q
$$
of ${\cal{I}}^q$ to the element
\[
  \left[ h(\omega) d_1^{\ast} s(\omega) - h(\omega) d_0^{\ast} s(\omega) \right] 
    \wedge d_1^{\ast} s(\omega_3 ) \wedge \cdots \wedge d_1^{\ast} s(\omega_q ) 
    \quad = \quad
    h(\omega) dh(\omega) \wedge \omega_3 \wedge \cdots \wedge \omega_q \pmod{ {\cal{B}}^{q-1} }
\]
which is a coboundary when $2$ is invertible ($d(h(\omega)^2) = 2 h(\omega)dh(\omega)$).

So, let us restrict our attention to the case $p = 2$ and show that the class $[\overline{\eta}|_{{\cal{I}}^q} ]$ is still zero. 
First, one can easily check that 
$\overline{\eta}: ( {\cal{H}}^1 )^{\otimes q} \to K^{q-1} / {\cal{B}}^{q-1}$,
and a fortiori $\overline{\eta}|_{ {\cal{I}}^{q}}: {\cal{I}}^q \to {\cal{Z}}^{q-1} / {\cal{B}}^{q-1} = {\cal{H}}^{q-1}$, is \emph{symmetric} in the sense that any element of the form
$$
\omega_1 \otimes \cdots \otimes \omega_{j-1} \otimes (\omega_j \otimes \omega_{j+1} + \omega_{j+1} \otimes \omega_{j} ) \otimes \omega_{j+2} \otimes \cdots \otimes \omega_q
$$
(when $1 = -1$, adding is subtracting) maps to zero under $\overline{\eta}$.
Therefore, one has a commutative diagram
$$
\xymatrix{
{\cal{I}}^{q} \ar@{>>}[r] \ar@{^{(}->}[d] & {\cal{I}}^q / {\cal{J}}^q \ar[r]^-{\widetilde{\eta}} \ar@{^{(}->}[d] & {\cal{H}}^{q-1} \ar[d] \ar@{^{(}->}[d] \\
({\cal{H}}^1)^{\otimes q} \ar@{>>}[r] & \Sym^q ( {\cal{H}}^1 ) \ar[r] & K^{q-1} / {\cal{B}}^{q-1}
}
$$
where the composite of the two horizontal arrows in the first row (resp. in the second row) is equal to $\overline{\eta}|_{{\cal{I}}^q}$  (resp. $\overline{\eta}$), and ${\cal{J}}^q$ denotes the (locally free) kernel of the projection $( {\cal{H}}^1)^{\otimes q} \to \Sym^q ( {\cal{H}}^1)$.

We get a notational advantage by taking the quotient by ${\cal{J}}^q$ : now ${\cal{I}}^q / {\cal{J}}^q$ is
generated by the images of local sections of the form
\begin{equation}
 \omega\otimes \omega \otimes \omega_3 \otimes \cdots \omega_q
 \label{eqn-s2.5}
\end{equation}
Such a local section is mapped under $\widetilde{\eta}$ onto
\begin{eqnarray}
 & & h(\omega) dh(\omega) \wedge d_0^{\ast} s (\omega_3 ) \wedge \cdots \wedge d_0^{\ast} s (\omega_q ) \pmod{ {\cal{B}}^{q-1}} \\
 & = & [ h(\omega) dh(\omega) \pmod{ {\cal{B}}^1 } ] \wedge \omega_3 \wedge \cdots \omega_q
 \label{eqn-s2.7}
\end{eqnarray}
We prove that $[\overline{\eta}|_{{\cal{I}}^q} ]$ is zero by finding a $0$-cochain $z$ with coefficients in $\underline{\Hom} ( {\cal{I}}^q , {\cal{H}}^{q-1} )$,
that is, a section of this sheaf over $U_0$, such that $\overline{\eta}|_{{\cal{I}}^q} = d_0^{\ast} z - d_1^{\ast} z$. 
As we know that $\overline{\eta}|_{{\cal{I}}^q}$ factors through $\widetilde{\eta}: {\cal{I}}^q / {\cal{J}}^q \to {\cal{H}}^{q-1}$, it suffices to find $\widetilde{z}: {\cal{I}}^q / {\cal{J}}^q \to {\cal{H}}^{q-1}$ such that
$$
\widetilde{\eta} = d_1^{\ast} \widetilde{z} - d_0^{\ast} \widetilde{z}.
$$
But from (\ref{eqn-s2.7}), one sees that
\begin{equation}
\widetilde{\eta} (\omega \otimes \omega \otimes \omega_3 \otimes \cdots \otimes \omega_q )
= \left( C^{-1} W^{\ast} dh(\omega) \right) \wedge \omega_3 \wedge \cdots \wedge \omega_q .
\label{eqn-s2.7bis}
\end{equation}
We denote here by $W$ the base-change of the absolute Frobenius endomorphism of $S$, so that the diagram
$$
\xymatrix{
X' \ar[r]^-{W} \ar[d] & X \ar[d] \\
S \ar[r]_-{\mathrm{frob}_S} & S
}
$$
is cartesian, by $W^{\ast}$ the pull-back morphism of differential forms
$$
\Omega^1_{X/S} \To W_{\ast} \Omega^{1}_{X'/S}
$$
and by $C^{-1}$ the (inverse) Cartier operation
$$
\Omega^1_{X'/S} \To {\cal{H}}^1 (F_{\ast} \Omega_{X/S}^{\bullet} )
$$
(\emph{cf.} \cite[\S 7]{Katz} and recall $p = 2$). Thus the last expression is the same as
$$
C^{-1} (d_1^{\ast} s\omega ) \wedge \omega_3 \wedge \cdots \wedge \omega_q
- C^{-1} (d_0^{\ast} s\omega ) \wedge \omega_3 \wedge \cdots \wedge \omega_q,
$$
and one is led to define over $U_0$
$$
\widetilde{z}: {\cal{I}}^q / {\cal{J}}^q \To {\cal{H}}^{q-1}
$$
so that it maps the local section (\ref{eqn-s2.5}) modulo ${\cal{J}}^q$ to
$$
C^{-1} ( s\omega ) \wedge \omega_3 \wedge \cdots \wedge \omega_q .
$$
As pointed out earlier, local sections of the form (\ref{eqn-s2.5}) generate ${\cal{I}}^q / {\cal{J}}^q$,
so such $\widetilde{z}$ is unique if exists at all. 
Now its existence can be shown locally: if one has a basis $e_1, \cdots, e_d$
of ${\cal{H}}^1$ over ${\cal{O}}_{X'}$, then the images of the sections
$$
\{ e_{j_1} \otimes \cdots \otimes e_{j_q} : 1 \le j_1 \le \cdots \le j_q \le d \,\, \mbox{ with at least one repetition } \,\, \}
$$
under ${\cal{I}}^q  \to {\cal{I}}^q / {\cal{J}}^q$ form a local basis of ${\cal{I}}^q / {\cal{J}}^q$, and then one can let $\widetilde{z}$ map the class of $e_{j_1} \otimes \cdots \otimes e_{j_q}$ to
$$
C^{-1} (s (\omega_{j_k} )) \wedge ( \mbox{ the rest } ),
$$
where $j_k$ is an index that repeats: If two or more indices repeat, whichever one is
chosen, the result is zero, and if an index repeats itself three or more times, it doesn't matter which consecutive terms are chosen, for $1 = -1$ and the sign doesn't matter.

Then one needs to show that any local section of the form (\ref{eqn-s2.5}) is mapped as
desired under $\widetilde{z}$ thus defined. 
One expresses the sections $\omega, \omega_3, \cdots , \omega_q$ as linear combinations of the $\{ e_i \}$ and one sees that it boils down to showing the linearity in each variable $\omega_3, \cdots, \omega_q$,
which is evident, as well as the linearity ``in the variable $\omega\otimes \omega$,'' which is less so.

Let $\omega = \alpha \xi + \beta \theta$, where $\alpha, \beta$ are sections of ${\cal{O}}_{X'}$ and $\xi, \theta$ sections of ${\cal{H}}^1$. Then one calculates
\begin{eqnarray*}
 C^{-1} W^{\ast} ( s(\alpha \xi + \beta \theta ) )
 & = & C^{-1} W^{\ast} ( F^{\ast} \alpha \cdot s(\xi ) + F^{\ast} \beta \cdot s(\theta) ) \\
 & = & C^{-1} ( \alpha^2 W^{\ast} s(\xi ) + \beta^2 W^{\ast} s(\theta ) ) \\
 & = & (F^{\ast} \alpha )^2 C^{-1} W^{\ast} s (\xi ) + (F^{\ast} \beta )^2 C^{-1} W^{\ast} s(\theta ),
\end{eqnarray*}
where $F^{\ast} : {\cal{O}}_{X'} \to F_{\ast} {\cal{O}}_X$ is the canonical pull-back morphism; here one uses the fact that $W\circ F$ is equal to the absolute second power Frobenius of $X$.

On the other hand, if one expands $\omega \otimes \omega$ as $\alpha^2 \xi\otimes \xi + \beta^2 \theta\otimes \theta + \alpha \beta (\xi \otimes \theta + \theta \otimes \xi)$,
then the last term is symmetric (i.e., lies in ${\cal{J}}^2$), and hence we get the same result this
way.

This can also be explained with the following diagram
$$
\xymatrix{
{\cal{Z}}^1 \subseteq F_{\ast} \Omega^1_{X/S} \ar[r]^-{F_{\ast} (W^{\ast})} &
        F_{\ast} W_{\ast} \Omega^1_{X'/S} \ar[r]^-{F_{\ast}W_{\ast} C^{-1}} & F_{\ast} W_{\ast} {\cal{H}}^1 (K) = (F_{X'})_{\ast} {\cal{H}}^1 (K) \\
{\cal{H}}^1 (K) \ar[u]^-{s}  & & 
}
$$
over $U_0$, where $F_{X'}$ denotes the absolute second power Frobenius of $X'$; it shows that
the map 
\[
    \omega\mapsto C^{-1} W^{\ast} s (\omega )
\]
is $2$-linear, while ``extracting'' $\omega$ out of $\omega\otimes \omega$ would be $2^{-1}$-linear; hence these nonlinearities cancel each other and the map $\omega\otimes \omega \mapsto C^{-1} W^{\ast} s(\omega )$ \emph{is} linear.

This shows that $\widetilde{z}$, hence the $0$-cochain $z$ which is obtained by composing $\widetilde{z}$ with the projection ${\cal{I}}^q \to {\cal{I}}^q / {\cal{J}}^q$, is well-defined and has the desired property. Therefore the class $[\overline{\eta}|_{{\cal{I}}^q}]$ is zero, and the only thing that contributes to the class (\ref{eqn-s2.4}) is the $q$th contraction of $\obs_1$. This ends the proof.
\end{proof}

By the construction (Theorem \ref{thm-s2.2.1}) and the calculation (Theorem \ref{thm-s2.3.1}), we immediately get:
\begin{corollary} 
 With the notations as in Theorem \ref{thm-s2.3.1}, suppose that $X' /S$ is liftable to $\widetilde{S}$.
Then for each integer $q$, the truncation $\tau_{[q-1, q]} F_{\ast} \Omega^{\bullet}_{X/S}$ of length $2$ is decomposable in the derived category $D(X', {\cal{O}}_{X'} )$.
\end{corollary}

\begin{proof}
This follows from \cite[3.5]{DeligneIllusie} (which identifies the obstruction to liftability to the decomposability of $\tau_{\leq 1}F_{X/S, \ast} \Omega^{\bullet}_{X/S}$)   and Theorems \ref{thm-s2.2.1} and \ref{thm-s2.3.1} (which relate the decomposability of $\tau_{[0,1]} F_{X/S, \ast} \Omega^{\bullet}_{X/S}$ with that of $\tau_{[q-1, q]} F_{X/S, \ast} \Omega^{\bullet}_{X/S}$).
\end{proof}

In particular, we extend the (special) case of Corollary~\ref{cor:two-step1} applied to the de Rham complex in characteristic $p>2$, even to the case of $p=2$.

\begin{corollary} \label{cor:delta-q-formula}
    Let $X$ be a smooth variety over a perfect field $k$.  Then we have the equality
  \[ 
    \delta^q_{F_{X/k,*}\Omega^\bullet_{X/k}} = (\id\otimes \delta^1_K)\circ \eta^q
  \]
  in $\Hom(\Omega^q_{X'/k}, \Omega^{q-1}_{X'/k}[2])$.
\end{corollary}

Finally, we answer the question of Katz:

\begin{corollary} \label{cor:katz-question}
  Let $S$ be a scheme of characteristic $p > 0$, $f : X \to S$ a smooth separated morphism of finite type, and $f' : X' \to S$ (resp. $F : X \to X'$) the
  base-change of $f$ by the Frobenius endomorphism of $S$ (resp. the relative Frobenius).
  Suppose $\widetilde{S}$ is a flat ${\mathbf{Z}} /p^2$-scheme whose reduction modulo $p$ yields $S$.
  Then the morphism in the conjugate spectral sequence
  $$
  d^{ij}_2 : R^i f'_{\ast} \Omega^j_{X'/S} \To R^{i+2} f'_{\ast} \Omega^{j-1}_{X'/S},
  $$
  where one identifies ${\cal{H}}^j (F_{\ast} \Omega_{X/S}^{\bullet} )$ with $\Omega^j_{X'/S}$ via the Cartier isomorphism, can be canonically regarded as the cup product with the additive inverse of the obstruction class (in $H^2 ( X', T_{X'/S} )$) to lifting $X'/S$ over $\widetilde{S}$.
\end{corollary}

\begin{proof}
We first remark that by \cite[3.9]{DeligneIllusie} it can be directly seen that the obstruction class to lifting does not depend on the choice of a flat $\mathbf{Z} /p^2$-lifting $\widetilde{S}$ of $S$. Then the corollary follows from \cite[3.5]{DeligneIllusie} and Theorems \ref{thm-s2.2.1} and \ref{thm-s2.3.1}.
\end{proof}

\appendix

\section{\texorpdfstring{$F$-split schemes of dimension $p+1$}{F-split schemes of dimension p+1}}

Let $k$ be a perfect field of characteristic $p>0$. As mentioned in the introduction, Drinfeld, Bhatt--Lurie, and Li--Mondal have obtained the following result.

\begin{theorem}[{\cite[\S 5.12.1]{Drinfeld}, \cite{BhattLurie}, \cite[Corollary~5.5]{LiMondal}}] \label{thm:drinfeld}
  Let $X$ be a smooth scheme over a perfect field $k$ of characteristic $p>0$. Suppose that $X$ is liftable to $W_2(k)$. Then the truncations 
  \[
    \tau_{[q-p+1, q]} F_{X/k,*} \Omega^\bullet_{X/k}
  \]
  are decomposable for all $q$.
\end{theorem}

Below, we employ this in order to show Kodaira--Akizuki--Nakano vanishing and Hodge--de Rham degeneration for $F$-split smooth projective schemes of dimension at most $p+1$. 

Recall \cite{MehtaRamanathan} that a $k$-scheme $X$ is $F$-split if the morphism $F_X^*\colon \cO_X\to F_{X,*} \cO_X$ is a split injection. Since $k$ is perfect, this is equivalent to the splitting of $F_{X/k}^*\colon \cO_{X'}\to F_{X/k,*} \cO_X$. It is well-known that every $F$-split scheme over $k$ admits a flat lifting to $W_2(k)$ \cite[\S 8.5]{IllusieFrobenius}.

If $X$ is $F$-split and if $L$ is a line bundle on $X$, then tensoring the split injection $\cO_X\to F_{X,*} \cO_X$ with $L$ and taking cohomology shows that for all $i$, $H^i(X, L)$ is a direct summand of $H^i(X, L\otimes F_* \cO_X)$. By the projection formula and the fact that Frobenius is affine this latter summand equals $H^i(X, F^* L) = H^i(X, L^p)$, and hence the Frobenius pull-back maps
\[ 
  F^* \colon H^i(X, L) \To H^i(X, L^p)
\]
are injective. Thus if $H^i(X, L^m) = 0$ for $m\gg 0$, then already $H^i(X, L) = 0$. Consequently, if $X$ is moreover smooth (or just Gorenstein) and projective, then $H^i(X, L^{-1}) = 0$ for $i<\dim X$ and $L$ ample, i.e.\ Kodaira vanishing holds on $X$. Similar reasoning with $L=\cO_X$ shows that 
\[
  F^* \colon H^i(X, \cO_X) \To H^i(X, \cO_X)
\]
is bijective for all $i\geq 0$.

\begin{theorem}[Kodaira--Akizuki--Nakano vanishing] \label{thm:KAN-van}
  Let $X$ be a smooth projective scheme over $k$ of dimension $d=p+1$. If $X$ is $F$-split, then Kodaira--Akizuki--Nakano vanishing holds for $X$, i.e.\ for every ample line bundle $L$, we have
  \[ 
    H^i(X, L^{-1}\otimes \Omega^j_{X/k}) = 0
    \quad\text{for}\quad
    i+j<d = p+1.
  \]
\end{theorem}

\begin{proof}
By Serre vanishing, the assertion holds for $L^{p^m}$ for $m\gg 0$. Therefore we may assume that it holds for $L^p$. Following \cite[Proof of Lemme~2.9]{DeligneIllusie}, we form the complex $K^\bullet = (L')^{-1}\otimes F_{X/k,*} \Omega^\bullet_{X/k}$ where $L'$ is the pull-back of $L$ to $X'$, and write the two spectral sequences
\begin{equation}
  {}_{I}E_1^{ij} = H^j(X', (L')^{-1}\otimes F_{X/k,*} \Omega^i_{X/k}) 
  \quad \Rightarrow \quad
  H^{i+j}(X', K^\bullet)
\end{equation}
and
\begin{equation} \label{eqn:conjss}
  {}_{II}E_2^{ij} = H^i(X', (L')^{-1}\otimes \Omega^j_{X'/k}) 
  \quad \Rightarrow \quad
  H^{i+j}(X', K^\bullet).
\end{equation}
Now the projection formula gives ${}_{I}E_1^{ij} = H^j(X, L^{-p}\otimes \Omega^i_{X/k})$, which vanishes for $i+j\leq p$ by assumption. Consequently the abutment $H^{r}(X', K^\bullet) = 0$ for $r\leq p$. 

We now investigate the second spectral sequence. Since $X$ is $F$-split, it lifts to $W_2(k)$. Theorem~\ref{thm:drinfeld} implies that the differentials on ${}_{II}E_r^{ij}$ are zero for $r\leq p$. For dimensional reasons, there are no nonzero differentials for $r>p+1$, and the only two nonzero differentials on ${}_{II}E_{p+1}^{ij}$ are
\[ 
  d_{p+1}^{0,p} \colon H^0(X', (L')^{-1} \otimes \Omega^p_{X'/k}) \To H^{p+1}(X', (L')^{-1})
\]
and
\[
  d_{p+1}^{0,p+1} \colon H^0(X', (L')^{-1} \otimes \omega_{X'/k}) \To H^{p+1}(X', (L')^{-1}\otimes \Omega^1_{X'/k}).
\]
We will show that $d_{p+1}^{0,p} = 0$, which will then imply that in \eqref{eqn:conjss} we have
\[
    {}_{II}E_2^{ij} = {}_{II}E_{p+1}^{ij} =  {}_{II} E_\infty^{ij}=0 
    \qquad \text{for $i+j\leq p$.}
\]
Note that $H^{p}(X', K^\bullet) = 0$ implies $0 = {}_{II}E_{p+2}^{0,p} = \ker(d_{p+1}^{0,p})$, i.e. $d_{p+1}^{0,p}$ is injective.

By Lemma~\ref{lem:commdiag} below applied to $E=L^{-1}$ and the map $F_{X^{-}/k}\colon X^- \to X$ where $X^- = (F_k)^{-1}(X)$ is the Frobenius \emph{untwist} of $X$, we have a commutative square induced by $F_{X/k}^*$:
\[ 
  \xymatrix{
    H^0(X', (L')^{-1} \otimes \Omega^p_{X'/k}) \ar[r]^-{d_{p+1}^{0,p}(L)} \ar[d] &  H^{p+1}(X', (L')^{-1}) \ar[d]^{F_{X'/k}^*}  \\
    0=H^0(X', L^{-p} \otimes \Omega^p_{X/k}) \ar[r]_-{d_{p+1}^{0,p}(L^p)} & H^{p+1}(X', L^{-p})
  }
\]
(in fact the left vertical map is zero). Here the bottom map is deduced similarly for (the pull-back to $X^{-}$ of) $L^p$ in place of $L$. The vertical right map is injective because $X$ is $F$-split, and the horizontal maps are injective by the previous paragraph. We conclude that $H^0(X', (L')^{-1} \otimes \Omega^p_{X'/k}) = 0$.
\end{proof}

In the proof above, as well as in the proof of Hodge--de Rham degeneration below, we need the following functoriality result.

\begin{lemma} \label{lem:commdiag}
  For a vector bundle $E$ on a smooth $k$-scheme $X$, write $K^\bullet(E)$ for the complex $E'\otimes F_{X/k,*} \Omega^\bullet_{X/k}$. Let $f\colon Y\to X$ be a map of smooth $k$-schemes, then $f$ induces a map of complexes $f^* K^\bullet(E)\to K^\bullet(f^* E)$ and hence a map of spectral sequences
  \[
    \xymatrix{
      {}_{II} E_2^{ij}(E) \ar@{=}[r] &  H^i(X', E'\otimes \Omega^j_{X'/k}) \ar@{}[r]|-{\displaystyle \Rightarrow} \ar[d] & H^{i+j}(X', K^\bullet(E)) \ar[d] \\
      {}_{II} E_2^{ij}(f^* E) \ar@{=}[r] &  H^i(Y', (f^* E)'\otimes \Omega^j_{Y'/k}) \ar@{}[r]|-{\displaystyle \Rightarrow} & H^{i+j}(Y', K^\bullet(f^* E)),
    } 
  \]
  where the maps $H^i(X', E'\otimes \Omega^j_{X'/k}) \to  H^i(Y', (f^* E)'\otimes \Omega^j_{Y'/k})$ are induced by the composition
  \[ 
    H^i(X', E'\otimes \Omega^j_{X'/k}) 
    \To 
    H^i(X', (f')^* (E')\otimes (f')^*\Omega^j_{X'/k})
    \To 
    H^i(X', (f')^* (E')\otimes \Omega^j_{Y'/k}).
  \] 
\end{lemma}

\begin{proof}
Obvious.
\end{proof}

\begin{theorem}[Hodge--de Rham degeneration]
  Let $X$ be a smooth projective scheme over $k$ of dimension $d=p+1$. If $X$ is $F$-split, then the Hodge to de Rham spectral sequence
  \[ 
    {}_I E_1^{ij} = H^j(X, \Omega^i_{X/k}) \quad \Rightarrow \quad H^{i+j}(X, \Omega^\bullet_{X/k})
  \]
  degenerates.
\end{theorem}

\begin{proof}
Since $X$ is proper, it is enough to show that the conjugate spectral sequence
\[ 
  {}_{II} E_2^{ij} = H^i(X', \Omega^j_{X'/k}) \quad \Rightarrow \quad H^{i+j}(X, \Omega^\bullet_{X/k})
\]
degenerates. Since $X$ is $F$-split, it lifts to $W_2(k)$, and then Theorem~\ref{thm:drinfeld} implies that the differentials on the page ${}_{II} E_r^{ij}$ of the conjugate spectral sequence are zero for $r\leq p$. 

Since $\dim X = p+1$, the only possibly nonzero differentials in this spectral sequence are therefore
\[ 
  d_{p+1}^{0,p} \colon H^0(X', \Omega^p_{X'/k}) \To H^{p+1}(X', \cO_{X'})
\]
and
\[ 
  d_{p+1}^{0,p+1} \colon H^0(X', \omega_{X'/k}) \To H^{p+1}(X', \Omega^1_{X'/k}).
\]
We will show that $d_{p+1}^{0,p} = 0$. Indeed, by functoriality of the above maps with respect to Frobenius (Lemma~\ref{lem:commdiag} with $E=\cO_X$ and the relative Frobenius $F_{X/k}$) gives a commutative square
\[ 
  \xymatrix{
    H^0(X', \Omega^p_{X'/k}) \ar[d]_{F_{X/k}^*} \ar[rr]^-{d_{p+1}^{0,p}} & & H^{p+1}(X', \cO_{X'}) \ar[d]^{F_{X/k}^*} \\
    H^0(X, \Omega^p_{X/k}) \ar[rr]_-{d_{p+1}^{0,p} \text{ for }X^-} & & H^{p+1}(X, \cO_{X})
  }
\]
where $X^- = (F^{-1}_k)^* X$ is again the Frobenius untwist of $X$. Since the Frobenius is zero on $\Omega^i$ for $i>0$, the left vertical map is zero. On the other hand, since $X$ is $F$-split, the right vertical map is an isomorphism. Therefore the top map $d_{p+1}^{0,p}$ is zero. 

Finally, we obtain the vanishing of $d_{p+1}^{0,p+1}$ by comparing dimensions and duality. Indeed, we have 
\begin{align*}
  \dim H^{p+2}(X, \Omega^\bullet_{X/k}) &= \dim H^{p}(X, \Omega^\bullet_{X/k}) & \text{(Poincar\'e duality)} \\
  &= \sum_{i+j=p} \dim H^{i}(X', \Omega^j_{X'/k}) & \text{(since $d_{p+1}^{0,p}=0$)} \\
  &= \sum_{i+j=p+2} \dim H^{i}(X', \Omega^j_{X'/k}). & \text{(Serre duality)}, 
\end{align*}
so $d_{p+1}^{0,p+1} = 0$.
\end{proof}

\begin{remark}[{See \cite[Theorem~5.4]{LiMondal} and \cite{BhattLurie}}]
In fact, the results of Drinfeld, Bhatt--Lurie, and Li--Mondal yield more than we have stated in Theorem~\ref{thm:drinfeld}. Namely, for $X$ smooth over $k$ and liftable to $W_2(k)$, there exists a decomposition in the derived category
\[ 
  F_{X/k,*} \Omega^\bullet_{X/k} \simeq \bigoplus_{i=0}^{p-1} K_i
\]
where $\h^j(K_i) = 0$ unless $i$ and $j$ are congruent modulo $p$. This implies that in the conjugate spectral sequence, as well as in the second spectral sequence used in the proof of Theorem~\ref{thm:KAN-van}), the only nonzero differentials may appear on pages $E_r$ where $r$ is congruent to one modulo $p$. We only used this with $r\leq p$, and it would be interesting to obtain new vanishing and degeneration theorems using this stronger fact.
\end{remark}

\bibliographystyle{amsalpha} 
\bibliography{bib}

\end{document}